\documentclass[11pt, a4paper]{article}
\usepackage{t1enc}
\usepackage[latin1]{inputenc}
\usepackage[english]{babel}
\usepackage{amsmath,amsthm}
\usepackage{amsfonts}
\usepackage{latexsym}
\usepackage[dvips]{graphicx}
\usepackage{float}
\usepackage[natural]{xcolor}
\usepackage{algorithm}
\usepackage{algorithmic}
\usepackage{enumerate}
\usepackage{multirow}
\usepackage{xcolor}
\usepackage[colorlinks,linkcolor=blue]{hyperref}
\usepackage{lineno}
\usepackage{setspace}
\usepackage{fullpage}

\usepackage{listings}
 \theoremstyle{plain}
\newtheorem{thm}{Theorem}[section]
\newtheorem{theorem}[thm]{Theorem}
\newtheorem{conjecture}[thm]{Conjecture}
\newtheorem{lemma}[thm]{Lemma}
\newtheorem{corollary}[thm]{Corollary}
\newtheorem{proposition}[thm]{Proposition}

\theoremstyle{definition}

\newtheorem{question}[thm]{Question}
\newtheorem{problem}[thm]{Problem}
\newtheorem{remark}[thm]{Remark}

\newtheorem{definition}[thm]{Definition}
\newtheorem{claim}[thm]{Claim}
\newtheorem{fact}[thm]{Fact}

\newtheorem{example}[thm]{Example}

\newtheorem{defn-thm}[thm]{Definition-Theorem}

\newcommand{\btheorem}{\begin{theorem}}
\newcommand{\etheorem}{\end{theorem}}
\newcommand{\bconjecture}{\begin{conjecture}}
\newcommand{\econjecture}{\end{conjecture}}
\newcommand{\bproposition}{\begin{proposition}}
\newcommand{\eproposition}{\end{proposition}}
\newcommand{\bdefinition}{\begin{definition}}
\newcommand{\edefinition}{\end{definition}}
\newcommand{\bcorollary}{\begin{corollary}}
\newcommand{\ecorollary}{\end{corollary}}
\newcommand{\bproof}{\begin{proof}}
\newcommand{\eproof}{\end{proof}}
\newcommand{\bclaim}{\begin{claim}}
\newcommand{\eclaim}{\end{claim}}
\newcommand{\bquestion}{\begin{question}}
\newcommand{\equestion}{\end{question}}
\newcommand{\bfact}{\begin{fact}}
\newcommand{\efact}{\end{fact}}
\newcommand{\bremark}{\begin{remark}}
\newcommand{\eremark}{\end{remark}}
\newcommand{\eexample}{\end{example}}
\newcommand{\bexample}{\begin{example}}

\newcommand{\elemma}{\end{lemma}}
\newcommand{\blemma}{\begin{lemma}}

\begin{document}

\title{Embedding clique-factors in graphs with low $\ell$-independence number}

\author{
Fan Chang\thanks{School of Mathematics, Shandong University, China. Email: {\tt fchang@mail.sdu.edu.cn, ghwang@sdu.edu.cn}. F.C. and G.W. were supported by Natural Science Foundation of China (11871311, 11631014) and Youth Interdisciplinary Innovation Group of Shandong University. }
 \and
Jie Han\thanks{Department of Mathematics, University of Rhode Island, RI, USA. Email: {\tt jie\_han@uri.edu.cn}.}
 \and
Jaehoon Kim\thanks{Department of Mathematical Sciences, KAIST, South Korea 34141. Email: {\tt jaehoon.kim@kaist.ac.kr}. J.K. was supported by the POSCO Science Fellowship of POSCO TJ Park Foundation and by the KAIX Challenge program of KAIST Advanced Institute for Science-X.}
 \and
Guanghui Wang\footnotemark[1]
 \and
 Donglei Yang\thanks{Data Science Institute, Shandong University, China, Email: {\tt dlyang@sdu.edu.cn}. D.Y. was supported by the China Postdoctoral Science Foundation (2021T140413) and Natural Science Foundation of China (12101365).}
}

\maketitle
\begin{abstract}

The following question was proposed by Nenadov and Pehova and reiterated by Knierim and Su: Given integers $\ell,r$ and $n$ with $n\in r\mathbb{N}$, is it true that every $n$-vertex graph $G$ with $\delta(G) \ge \max \{ \frac{1}{2},\frac{r - \ell}{r} \}n + o(n)$ and $\alpha_{\ell}(G) = o(n) $ contains a $K_{r}$-factor?
We give a negative answer for the case when $\ell\ge \frac{3r}{4}$ by giving a family of constructions using the so-called cover thresholds
and show that the minimum degree condition given by our construction is asymptotically best possible.
That is, for all integers $r,\ell$ with $r > \ell \ge \frac{3}{4}r$ and $\mu >0$, there exist $\alpha > 0$ and $N$ such that for every $n\in r\mathbb{N}$ with $n>N$, every $n$-vertex graph $G$ with $\delta(G) \ge \left( \frac{1}{2-\varrho_{\ell}(r-1)} + \mu \right)n $ and $\alpha_{\ell}(G) \le \alpha n$ contains a $K_{r}$-factor.
Here $\varrho_{\ell}(r-1)$ is the Ramsey--Tur\'an density for $K_{r-1}$ under the $\ell$-independence number condition.

\end{abstract}

\section{Introduction}

In this paper we study the following problem proposed by Nenadov and Pehova~\cite{Nenadov2018}.

\begin{problem}\label{prob1.3}
Let $r,\ell \in \mathbb{N}$ with $r \ge \ell \ge 2$ and $G$ be an $n$-vertex graph with $n\in r\mathbb{N}$ and $\alpha_{\ell}(G) = o(n)$. What is the best possible minimum degree condition on $G$ that guarantees a $K_{r}$-factor?
\end{problem}

Here given a graph $G$, $\alpha_{\ell}(G)$ denotes the $\ell$-\emph{independence number} of $G$, the size of a largest $K_{\ell}$-free induced subgraph of $G$ (the independence number of $G$ corresponds to $\alpha_{2}(G)$).

\subsection{Hajnal--Szemer\'edi Theorem}

Given graphs $H$ and $G$, a collection of vertex-disjoint copies of $H$ in $G$ is called an $H$-\emph{tiling}. A \emph{perfect $H$-tiling} of $G$, or an \emph{$H$-factor}, is an $H$-tiling which covers all the vertices of $G$.

Determining sufficient conditions for the existence of an $H$-factor is one of the fundamental lines of research in extremal graph theory -- one prominent reason is due to a result of Hell and Kirkpatrick~\cite{1983Hell} which shows that the decision problem for $H$-factors is NP-complete, given that $H$ has a connected component of size at least 3.
The first result in this direction is by Dirac~\cite{Dirac1952}, who showed that an $n$-vertex graph $G$ with minimum degree at least $\frac{n}{2}$ contains a Hamilton cycle, in particular if $n$ is even then $G$ has a perfect matching.
Then a celebrated result of Hajnal and Szemer\'edi~\cite{hajnal1970proof} gives that $\delta(G) \ge \frac{(r - 1)}{r}n$ is sufficient for the existence of a $K_{r}$-factor (the $r=3$ case was previously obtained by Corr\'adi and Hajnal~\cite{K1963On}).
The complete $r$-partite graph with each part of size $\frac{n}{r}$ shows the tightness of the minimum degree condition.
Determining the best possible bound on $\delta(G)$ for an arbitrary $H$ has been a major problem and intrigued many excellent works (see e.g.~\cite{alon1996h,komlos2000tiling,komlos2001proof,shokoufandeh2003proof} and the survey~\cite{D2009Surveys}) until it was finally settled by K\"uhn and Osthus~\cite{D2009The}.

There have been various generalizations of Hajnal--Szemer\'edi Theorem in the setting of partite graphs~\cite{keevash2015multipartite,magyar2002tripartite,martin2008quadripartite}, directed graphs~\cite{Treglown2015}, and hypergraphs~\cite{rodl2009perfect}. For each of these results, the extremal graphs all have one or more large independent sets.
A natural follow-up question is then how much the bound on the minimum degree can be weakened if we forbid large independent sets?
The problem was first studied by Balogh, Molla and Sharifzadeh~\cite{Balogh2016Triangle} who solved the problem for triangles -- they proved that for $n$-vertex graphs with sublinear independence number, the minimum degree condition forcing a $K_3$-factor can be weakened from $\frac{2n}{3}$ to $\frac{n}{2}+o(n)$.
Nenadov and Pehova~\cite{Nenadov2018} generalized this problem to $\ell$-independence numbers for $\ell \ge 2$, that is, they proposed Problem~\ref{prob1.3}.
They also proved a general upper bound for the minimum degree conditions, which is asymptotically optimal when $\ell=k-1$.
Recently an excellent work of Knierim and Su~\cite{knierim2019kr} resolved the case $r\ge 4$ and $\ell=2$ by showing that $\delta(G) \ge (\frac{r - 2}{r}+o(1))n$ is sufficient.
They reiterated Problem~\ref{prob1.3} in their paper and proposed a minimum degree condition as follows.



\begin{problem}\label{prob1.5}\cite{knierim2019kr,Nenadov2018}
Is it true that for every $r,\ell\in\mathbb{N}$ with $r\ge \ell$ and $\mu>0$ there is a constant $\alpha$ and $n_0 \in \mathbb{N}$ such that every graph $G$ on $n\ge n_0$ vertices where $r$ divides $n$ with $\delta(G) \ge \max \{ \frac{1}{2}+\mu,\frac{r - \ell}{r}+\mu \}n$ and $\alpha_{\ell}(G) \le \alpha n $ has a $K_{r}$-factor?
\end{problem}

In this paper we determine the asymptotically optimal minimum degree condition for $\ell \ge \frac{3}{4}r$, which solves Problem~\ref{prob1.3} for this range, and actually provides a negative answer to Problem~\ref{prob1.5}.
\btheorem\label{main thm}
Let $r,\ell \in \mathbb{N}$ such that $r > \ell \ge \frac{3}{4}r$. For any $\mu >0$, there exists an $\alpha > 0$ such that for all sufficiently large $n\in r\mathbb{N}$, every $n$-vertex graph $G$ with $\delta(G) \ge \left( \frac{1}{2-\varrho_{\ell}(r-1)} + \mu \right)n $ and $\alpha_{\ell}(G) \le \alpha n$  contains a $K_{r}$-factor.
Moreover, the minimum degree condition is asymptotically best possible.
\etheorem

The constant $\varrho_{\ell}(r-1)$ will be formally defined in the next subsection.

\subsection{Ramsey--Tur\'an Theory}\label{sec1.2}

Initiated in 1969 by S\'os, and later generalized by Erd\H{o}s, Hajnal, S\'os and Szemer\'edi~\cite{1983More}, \emph{Ramsey--Tur\'an theory} studies the maximum number of edges in an $n$-vertex $K_r$-free graph with sublinear independence number. Erd\H{o}s and S\'os~\cite{1970More} solved this problem when $r$ is an odd integer. The case when $r$ is even is more intricate. Let $\varrho_{\ell}(r)$ be the \emph{Ramsey--Tur\'an density}: $\varrho_{\ell}(r):=\lim\limits_{\varepsilon\to0}\lim\limits_{n\to\infty}\frac{\textbf{RT}_{\ell}(n, K_r, \varepsilon n)}{\binom{n}{2}}$, where $\textbf{RT}_{\ell}(n, K_r, \varepsilon n)$ denotes the maximum number of edges of an $n$-vertex $K_r$-free graph $G$ with $\alpha_{\ell}(G) \le \varepsilon n$. Szemer\'edi~\cite{E1972On} first showed that $\textbf{RT}(2, 4)\le\frac{1}{4}$. This turned out to be sharp as four years later Bollob\'as and Erd\H{o}s~\cite{1976On} provided a matching lower bound  using an ingenious geometric construction. There are some recent exciting developments in this area~\cite{BaloghRT2012,BaloghRT2013,Liu2021}. For further information on Ramsey--Tur\'an theory the reader is referred to the comprehensive survey~\cite{2001Ramsey} by Simonovits and S\'os.

For our purpose, we define the minimum-degree type analogue for the Ramsey--Tur\'an density as follows.
Here for $\alpha > 0$, we denote by $\textbf{RT}^*_{\ell}(n,K_{r},\alpha n)$ the minimum integer $\delta$ such that every $n$-vertex graph $G$ with $\delta(G) \ge \delta$ and $\alpha_{\ell}(G)\le \alpha n$ contains a copy of $K_{r}$, and thus let \[\varrho^*_\ell(r):=\lim\limits_{\alpha\to 0}\liminf\limits_{n\to \infty}\frac{\textbf{RT}^*_{\ell}(n,K_{r},\alpha n)}{n}.\] We write $\textbf{RT}^*_{\ell}(n,K_{r},o(n))=\varrho^*_\ell(r)n+o(n)$.
Moreover, the following equality always holds.

\begin{proposition}\label{RTprop}
Let $r,\ell\in\mathbb{N}$ with $r\ge\ell\ge2$. Then $\varrho^*_\ell(r)=\varrho_\ell(r)$.
\end{proposition}

Note that the following result due to Balogh and Lenz~\cite{BaloghRT2013} shows that $\varrho_\ell(r)>0$ for all $r\ge\ell+2$.

\btheorem[Balogh and Lenz~\cite{BaloghRT2013}]\label{RTthm}
For $\ell\ge2$ and $\ell+2\le r\le 2\ell$, let $u=\lceil\frac{\ell}{2}\rceil$. Then
$$
\varrho_\ell(r)\ge\left(1-\frac{1}{r-\ell}\right)2^{-u^2}.
$$
\etheorem

\subsection{Cover Threshold}
For our problem, the ``cover'' property plays an essential role. To have a $K_r$-factor in $G$, we clearly have to be above the \emph{cover threshold}: every vertex $v\in V(G)$ must be covered by a copy of $K_r$ in $G$. The definition of $\varrho^*_{\ell}(r-1)$ and Proposition~\ref{def} actually tell us that every $v\in V(G)$ is covered by many copies of $K_r$ in $G$. The cover threshold has been first discussed in~\cite{Han2017} and appeared in a few different contexts~\cite{chang2020factors,Sun2021quasirandom}.

\begin{figure}[htb]
\center{\includegraphics[width=6.6cm] {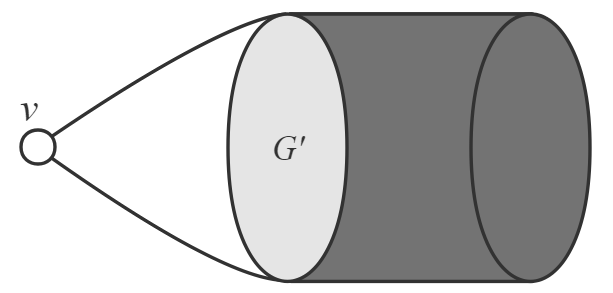}}
\caption{An extremal graph where $G'$ is $K_{r-1}$-free.}
\label{figure1}
\end{figure}

Now we give a construction that shows the optimality of Theorem~\ref{main thm} (the moreover part of the theorem). For any $x\in(0,1]$, let $G$ be an $n$-vertex graph defined by firstly fixing a vertex $v$ such that $N(v)=xn$ and $G':=G[N(v)]$ is a $K_{r-1}$-free subgraph with $\delta(G') = \varrho^*_{\ell}(r-1)xn-o(n)$ and $\alpha_{\ell}(G') = o(n)$, and then adding a clique of size $(n-xn-1)$ that is complete to $N(v)$ (see Figure~\ref{figure1}). Here there exists no copy of $K_r$ covering $v$ and thus $G$ contains no $K_r$-factor; moreover, by choosing $x=\frac{1}{2-\varrho^*_{\ell}(r-1)}$ we obtain $\delta(G)=\frac{1}{2-\varrho^*_{\ell}(r-1)}n-o(n)$. This construction combined with Theorem~\ref{RTthm} gives a negative answer to Problem~\ref{prob1.5}.

Note that the definition of $\varrho^*_{\ell}(r-1)$ implies the following proposition, which will be used in later proofs.

\begin{proposition}\label{def}
Given $r,\ell \in \mathbb{N}$ and a constant $\mu >0$, there exists $\alpha >0$ such that for all sufficiently large $n$ the following holds. Let $G$ be an $n$-vertex graph with $\delta(G) \ge \left( \frac{1}{2-\varrho^*_{\ell}(r-1)} + \mu \right)n$ and $\alpha_{\ell}(G)\le\alpha n$. If $W$ is a subset of $V(G)$ with $|W|\le\frac{\mu}{2}n$, then for each vertex $u\in V(G)\setminus W$, $G[V(G)\setminus W]$ contains a copy of $K_r$ covering $u$.
\end{proposition}
\bproof
Let $G_1:=G[V(G)\setminus W]$. It suffices to show that for each vertex $u\in G_1$, there is a copy of $K_{r-1}$ in $N_{G_1}(u)$. Note that for every vertex $u$ in $G_1$, we have $|N_{G_1}(u)|\ge \delta(G_1)\ge\delta(G)-|W|\ge\left( \frac{1}{2-\varrho^*_{\ell}(r-1)} + \frac{\mu}{2} \right)n$. Given any vertex $v\in N_{G_1}(u)$, letting $d_{G_1}(u,v):=|N_{G_1}(u)\cap N_{G_1}(v)|$, we have
\[
  \begin{split}
d_{G_1}(u,v)-\left(\varrho^*_{\ell}(r-1)+\frac{\mu}{4} \right)d_{G_1}(u)&\ge d_{G_1}(u)+d_{G_1}(v)-n-\left(\varrho^*_{\ell}(r-1)+\frac{\mu}{4} \right)d_{G_1}(u)\\
&\ge\left(2-\varrho^*_{\ell}(r-1)-\frac{\mu}{4} \right)\delta(G_1)-n>\frac{\mu}{8}n >0.
 \end{split}
 \]
Thus $\delta(G[N_{G_1}(u)])>(\varrho^*_{\ell}(r-1)+\frac{\mu}{4})|N_{G_1}(u)|$. Therefore by the definition of $\varrho^*_{\ell}(r-1)$, $G[N_{G_1}(u)]$ contains a copy of $K_{r-1}$ since $\alpha$ is sufficiently small (comparing with $\mu$), which together with $u$ yields a copy of $K_r$ in $G_1$.
\eproof

\subsection{Comparison with related results}

Another way to consider our problem is to formulate the question as adding or requiring randomness in host graphs.
In this way, comparison can be drawn to other related problems, namely, the randomly perturbed models and the locally dense model.
The randomly perturbed graphs, firstly studied by Bohman, Frieze and Martin~\cite{Randomperturbed2003}, consider adding random edges to a deterministic graph with a minimum degree condition and have been studied extensively in recent years.
The locally dense graphs can be traced back to Erd\H{o}s and S\'os~\cite{ES1982} which proposed a study of graphs with uniform edge distribution.

Let us draw a comparison with our problem where for instance we consider $\alpha_2(G)=o(n)$.
For the randomly perturbed model this is equivalent to adding linearly many random edges which already brings down the minimum degree requirement for many properties significantly, most notably, the Hamiltonicity (the minimum degree requirement drops from $\frac{n}{2}$ to $o(n)$).
In comparison, adding the sublinear independence number condition does not reduce the requirement for the Hamiltonicity (e.g., consider two disjoint cliques of equal sizes).
Further comparison on clique-factor problems can be done to the work of~\cite{HMT2021}.
On the other hand, as observed by Reiher and Schacht~\cite{reiher_schacht_2019}, the locally dense condition is quite handy in the graph setting, which, e.g., allows one to draw an arbitrarily large constant-sized clique in a linear-sized set of vertices. Based on this, they showed that locally dense graphs with positive density and minimum degree at least $\frac{n}{2}+o(n)$ contain $K_r$-factors for every fixed $r$ (given that $r\mid n$).
Therefore, the sublinear $\ell$-independence number condition is probably the weakest condition among all kinds of ``randomness'' assumptions.

\subsection{Proof strategy}

Our proof uses the framework of the absorption method, pioneered by the work of R\"odl, Ruci\'nski and Szemer\'edi~\cite{rodl2009perfect} on perfect matchings in hypergraphs (though similar ideas appear implicitly in previous works, see e.g.~\cite{Krivelevich_TF}). 
This method is particularly efficient in constructing spanning subgraphs, before its appearance the only known systematic way is the blow-up lemma due to Koml\'os, S\'ark\"ozy and Szemer\'edi~\cite{KSS_blowup}.
A typical key step in the absorption method for $H$-factor problem is to show that for every set of $h:=|V(H)|$ vertices, the host graph $G$ contains $\Omega(n^b)$ $b$-vertex absorbers (to be defined shortly).
However,  as pointed out in~\cite{Balogh2016Triangle}, in our setting this is usually impossible because when we construct the absorbers using the independence number condition, it does not give such a strong counting.
Instead, a much weaker notion is used in this paper, that is, we aim to show that \emph{for (almost) every set of $h$ vertices, the host graph $G$ contains $\Omega(n)$ vertex-disjoint absorbers}.
Note that this weak notion of absorbers have been successfully used in our setting~\cite{Nenadov2018, knierim2019kr} and the randomly perturbed setting~\cite{chang2020factors}.

\subsubsection{The absorption method}
Given the above weakening on the notion of absorbers, the general proof strategy follows the original work-flow of the absorption method:
the main tasks are to ($i$) build an absorbing set $A$ and ($ii$) find an $H$-tiling in $G[V(G)\setminus A]$ which covers almost all vertices.

To formulate $(i)$ we first give the following definitions.

\bdefinition Let $H$ be a graph with $h$ vertices and $G$ be a graph with $n$ vertices.
\begin{enumerate}
	\item We say that a subset $A \subseteq V(G)$ is a $\xi$\emph{-absorbing set} for some $ \xi > 0$ if for every subset $R \subseteq V(G)\setminus A$ with $|R| \le \xi n$ and $|A\cup R|\in h\mathbb{N}$, $G[A \cup R]$ contains an $H$-factor.
	\item Given a subset $S \subseteq V(G)$ of size $h$ and an integer $t$, we say that a subset $A_{S} \subseteq V(G)\setminus S$ is an $(S,t)$\emph{-absorber} if $|A_{S}| = ht$ and both $G[A_{S}]$ and $G[A_{S} \cup S]$ contain an $H$-factor.
\end{enumerate}
\edefinition

One has to show that one can still construct an absorbing set, given only the weak notion of absorber counting as discussed above.
This is now known possible via the clever idea ``bipartite template'' first used by Montgomery~\cite{Montgomery}.
Such an approach is summarized in our setting as the following result by Nenadov and Pehova.

\blemma[Absorbing lemma~\cite{Nenadov2018}]\label{Nenadov lem}
Let $H$ be a graph with $h$ vertices and let $\gamma > 0$ and $t \in \mathbb{N}$ be constants.  Then there exist $\xi = \xi(h,t,\gamma)$ and $n_0\in \mathbb{N}$ such that the following statement holds. Suppose that $G$ is a graph with $n \ge n_{0}$ vertices such that every $S \in \tbinom{V(G)}{h} $ has a family of at least $\gamma n$ vertex-disjoint $(S,t)$-absorbers. Then $G$ contains a $\xi$-absorbing set of size at most $\gamma n$.
\elemma

Hence, our main task is to find linearly many vertex-disjoint absorbers for every $h$-subset.
In order to achieve this purpose, we shall use the \emph{lattice-based absorbing method} developed recently by the second author~\cite{Han2017Decision}.
The reachability argument was introduced by Lo and Markstr\"om~\cite{2015Lo}, and we here state a slightly different version reflecting our weak notion of absorbers.
Let $G$ be a graph of $n$ vertices and $H$ be a graph of $h$ vertices. For any two vertices $u, v \in V(G)$ and an integer $t \ge 1$, a set $S$ of size at most $ht-1$ is called an \emph{$H$-reachable set} for $\{u, v\}$ if both $G[\{u\} \cup S]$ and $G[\{v\} \cup S]$ have $H$-factors. For $ t \ge 1$ and $\beta > 0$, we say that two vertices $u$ and $v$ are \emph{$(H,\beta,t)$-reachable} (in $G$) if there are $\beta n$ vertex-disjoint $H$-reachable sets $S$ in $G$. Moreover, we say that a vertex set $U \subseteq V(G)$ is \emph{$(H,\beta,t)$-closed} if any two vertices in $U$ are $(H,\beta,t)$-reachable in $G$.

To verify $(ii)$, our main tools are the regularity lemma, a result of Koml\'os~\cite{komlos2000tiling} , and embedding techniques in the setting of low $\ell$-independence number (see Lemma~\ref{tiling lem}).

\subsubsection{Main technical challenge -- Lemma~\ref{reachable lem}}
Under the reachability framework, a key step is to show that after excluding a small set $B$ of vertices, every vertex $v$ in $G':=G[V(G)\setminus B]$ is still $(K_r,\beta_1,1)$-reachable (in $G'$) to many other vertices $u$ in $G'$ (see Lemma~\ref{reachable lem}).
In order to prove Lemma~\ref{reachable lem}, we use the regularity method which helps us to capture the rough structure of the graph.
In an $\varepsilon$-regular partition, naively typical pairs of vertices from the same cluster should behave ``alike'' each other.
Thus, to establish Lemma~\ref{reachable lem} we choose $u$ as in the cluster in which $v$ is and use the regularity method to obtain a ``clean'' structure of $G$ so that we can establish a small induced subgraph of $G$ in $N(u)\cap N(v)$ with a high minimum degree and then apply the definition of $\rho(\ell, r-1)$ (as in the proof of Proposition~\ref{def}).
There are some other technicalities in the proof, e.g., to control the regularity, we have to work with a constant number of clusters instead of all of them, and thus we need to show that a minimum-degree type estimate is inherited by a random subset of clusters by standard concentration method (which unfortunately incurs some lengthy calculations).





\medskip\noindent\textbf{Notation.} Throughout the paper we follow the standard graph-theoretic notation~\cite{Diestel2017}. For a graph $G=(V,E)$, we let $v(G)=|V|$ and $e(G)=|E|$. For $U\subseteq V$, $G[U]$ denotes the induced graph of $G$ on $U$. The notation $G-U$ is used to denote the induced graph after
removing $U$, that is $G-U:=G[V\setminus U]$. For two subsets $A,B\subseteq V(G)$, we use $e(A,B)$ to denote the number of edges joining $A$ and $B$. Given a vertex $v\in v(G)$ and $X\subseteq V(G)$, denote by $N_X(v)$ the set of neighbors of $v$ in $X$ and let $d_X(v):=|N_X(v)|$. In particular, we write $N_{G}(v)$ for the set of neighbors of $v$ in $G$. We omit the index $G$ if the graph is clear from the context. Given a set $V$ and an integer $k$, we write $\binom{V}{k}$ for the family of all $k$-subsets of $V$. For any integers $a\le b$, let $[a,b]:=\{i\in\mathbb{Z}:a\le i\le b\}$.

When we write $\alpha \ll \beta\ll \gamma$, we always mean that $\alpha, \beta, \gamma$ are constants in $(0,1)$, and $\beta\ll \gamma$ means that there exists $\beta_0=\beta_0(\gamma)$ such that the subsequent arguments hold for all $0<\beta\le \beta_0$.
Hierarchies of other lengths are defined analogously.

\medskip\noindent\textbf{Organization.} The rest of the paper is organized as follows. In the next section, we shall introduce an absorbing lemma (Lemma~\ref{absorbing lem}) and an almost perfect tiling lemma (Lemma~\ref{tiling lem}) to prove Theorem~\ref{main thm}. In Section~\ref{sec3}, we shall prove Lemma~\ref{tiling lem}. In Section~\ref{sec4}, we will briefly present some necessary results and tools to introduce the lattice-based absorbing method and prove Lemma~\ref{absorbing lem}. We shall prove Lemma~\ref{reachable lem} in Section~\ref{sec5} and Lemma~\ref{merge lem} in Section~\ref{sec6}. In Section~\ref{sec8}, we shall give a proof of Proposition~\ref{RTprop}, which is relatively independent. In Section~\ref{sec7} we shall propose some questions.

\section{Proof of Theorem~\ref{main thm} and Preliminaries}

\subsection{Proof of Theorem~\ref{main thm}}
Following the popular absorption method, our main work is to establish an absorbing lemma (Lemma~\ref{absorbing lem}) and an almost-perfect tiling lemma (Lemma~\ref{tiling lem}), which together yield Theorem~\ref{main thm}.
Below we state these two lemmas, whose roles hopefully have been motivated well in Section~\ref{sec1.2}.

\blemma[Absorbing Lemma]\label{absorbing lem}
Given positive integers $r,\ell$ with $r>\ell>\frac{2}{3}r$ and positive constants $\mu,\gamma$ with $\gamma\le\frac{\mu}{2}$, there exist $\alpha>0$, $\xi>0$ and $n_0\in \mathbb{N}$ such that the following holds for every $n\ge n_0$. Let $G$ be an $n$-vertex graph with $\delta(G) \ge \left( \frac{1}{2-\varrho^*_{\ell}(r-1)} + \mu \right)n $ and $\alpha_{\ell}(G) \le \alpha n$. Then $G$ contains a $\xi$-absorbing set $A$ of size at most $\gamma n$.

\elemma

Lemma~\ref{tiling lem} provides a $K_r$-tiling that covers almost all vertices, whose proof will be given in Section~\ref{sec3}.

\blemma[Almost perfect tiling]\label{tiling lem}
Given positive integers $r,\ell$ such that $r > \ell \ge \frac{3}{4}r$ and positive constants $\mu,\delta$, there exist $\alpha> 0$ and $n_{0} \in \mathbb{N}$ such that every graph $G$ on $n \ge n_{0}$ vertices with $\delta(G) \ge \left( \frac{1}{2} + \mu \right)n $ and $\alpha_{\ell}(G) \le \alpha n$ contains a $K_{r}$-tiling covering all but at most $\delta n $ vertices in $G$.
\elemma

Now we are ready to prove Theorem~\ref{main thm} using Lemma~\ref{absorbing lem} and Lemma~\ref{tiling lem}.

\bproof[Proof of Theorem~\ref{main thm}]

Given any positive integers $\ell,r$ with $r>\ell\ge\frac{3r}{4}$ and a constant $\mu>0$.  Choose $\frac{1}{n}<\alpha\ll\delta\ll\xi\ll\gamma\ll\mu$. Let $G$ be an $n$-vertex graph with $\delta(G)\ge \left( \frac{1}{2-\varrho^*_{\ell}(r-1)} + \mu \right)n $ and $\alpha_{\ell}(G) \le \alpha n$.

By Lemma~\ref{absorbing lem} with $\gamma\le\frac{\mu}{2}$, we find a $\xi$-absorbing set $A \subseteq V(G)$ of size at most $\gamma n$. Let $G_1:=G-A$. Then we have
\[
\delta(G_1) \ge \left( \frac{1}{2-\varrho^*_{\ell}(r-1)} + \mu \right)n-\gamma n\ge \left( \frac{1}{2-\varrho^*_{\ell}(r-1)} + \frac{\mu}{2} \right)n.
\]
Therefore applying Lemma~\ref{tiling lem} on $G_1$ with $\delta$, we obtain a $K_{r}$-tiling $\mathcal{M}$ that covers all but a set $R$ of at most $\delta n$ vertices in $G_1$. Since $\delta\ll\xi$, the absorbing property of $A$ implies that $G[A \cup R]$ contains a $K_r$-factor $\mathcal{R}$, which together with $\mathcal{M}$ forms a $K_r$-factor in $G$.
\eproof

\subsection{Regularity}
An important component of our proof is the Szemer\'edi's Regularity Lemma.
We begin with some notation. Given a graph $G$ and a pair $(A, B)$ of vertex-disjoint subsets in $V(G)$, the \emph{density} of $(X,Y)$ is defined as
\[
d(X,Y) = \frac{e(X,Y)}{|X||Y|}.
\]
Given a constant $\varepsilon >0$, we say that $(X, Y)$ is an $\varepsilon $-\emph{regular pair} in $G$ (or $(X, Y)$ is $\varepsilon $-\emph{regular}) if for all $X' \subseteq X$, $Y' \subseteq Y$ with  $|X'| \ge \varepsilon |X|$ and $|Y'| \ge \varepsilon |Y|$ we have
\[
|d(X',Y') - d(X,Y)|  \le  \varepsilon.
\]

Moreover, a pair $(X, Y)$ is called $(\varepsilon,\delta)$-$super$-$regular$ for some $\delta>0$ if $(X, Y)$ is $\varepsilon $-regular and $d_{Y}(x) > \delta|Y|$ for all $x \in X$ and $d_{X}(y) > \delta|X|$ for all $y \in Y$. The following fact is an easy consequence of the definition of regularity.

\bfact\label{slicing lem}
Given constants $\eta >\varepsilon>0$, and a bipartite graph $G =(X \cup Y, E)$, if $(X,Y)$ is $\varepsilon $-regular, then for all $X_{1} \subseteq X$ and $Y_{1} \subseteq Y$ with $|X_{1}| \ge \eta|X|$ and $|Y_{1}| \ge \eta|Y|$, we have that $(X_{1}, Y_{1})$ is $\varepsilon '$-regular in $G$ for any $\varepsilon ' \ge {\rm max}\{ \frac{\varepsilon}{\eta},2\varepsilon  \}$. 
\efact

We now state the degree form of the regularity lemma (see~\cite[Theorem 1.10]{komlos1996}).

\blemma[Degree form of the Regularity Lemma~\cite{komlos1996}]\label{reg}

For every $\varepsilon > 0$ there is an $N = N(\varepsilon )$ such that the following holds for any real number $d \in [0,1]$ and $n\in \mathbb{N}$. Let $G=(V,E)$ be a graph with $n$ vertices. Then there exists a partition $\mathcal{P}=\{V_{0},\ldots, V_{k}\} $ of $V$ and a spanning subgraph $G' \subseteq G$ with the following properties:
\begin{enumerate}
  \item [$(a)$]\label{a} $ \frac{1}{\varepsilon}\le k \le N $;
  \item [$(b)$] $|V_{i}| \le \varepsilon  n $ for $0 \le i \le k$ and $|V_{1}|=|V_{2}|=\cdots=|V_{k}| =m$ for some $m\in \mathbb{N}$;
  \item [$(c)$] $d_{G'}(v) > d_{G}(v) - (d + \varepsilon )n$ for every $v \in V(G)$;
  \item [$(d)$] every $V_{i}$ is an independent set in $G'$;
  \item [$(e)$] each pair $(V_{i},V_{j})$, $1 \le i < j \le k$ is $\varepsilon $-regular in $G'$ with density 0 or at least $d$.
\end{enumerate}
\elemma

A widely-used auxiliary graph accompanied with the regular partition is the reduced graph.
The \emph{$d$-reduced graph} $R_d$ of $\mathcal{P}$ is a graph defined on the vertex set $\{V_1,\ldots,V_k\}$ such that $V_i$ is connected to $V_j$ by an edge if $(V_i,V_j)$ has density at least $d$ in $G'$. To ease the notation, we use $d_R(V_i)$ to denote the degree of $V_i$ in $R_d$ for each $i\in[k]$. In our proof, we also need a weighted version as follows. A \emph{weighted reduced graph} of $\mathcal{P}$, denoted by $R^*_d$, is obtained from $R_d$ by adding a weight function $d: E(R_d)\rightarrow [0,1]$ which assigns to every edge $V_iV_j$ the density $d(V_i,V_j)$ in $G'$. To simplify the notation, for every edge $V_iV_j$, we write $d_{ij}$ for the weight $d(V_i,V_j)$  and let $d_{R^*}(V_i) := \sum_{V_j \sim V_i} d_{ij}$. 
\begin{fact}\label{min degree}
Given positive constants $d,\mu,\varepsilon$ and $c$, fix an $n$-vertex graph $G=(V,E)$ with $\delta(G) \ge ( c + \mu)n$ and let $G'$ and $\mathcal{P}$ be obtained by Lemma~\ref{reg}, and $R^*_d$ be given as above. Then for every $V_i\in V(R^*_d)$ we have
$d_{R^*}(V_i)\ge (c + \mu -2\varepsilon  -d)k.$
\end{fact}
\bproof
Note that $|V_0|\le \varepsilon n$ and $|V_i|=m$ for each $V_i\in V(R^*_d)$. Thus we have that
\[
\sum_{V_j \sim V_i} d_{ij}|V_i||V_j|=e_{G'}(V_{i},\cup_{j \ne i}V_{j})\ge (\delta(G')-|V_0|)|V_i|\ge\left(c + \mu -2\varepsilon  -d \right)nm,
\]
which implies that
\[
d_{R^*}(V_i)=\sum_{V_j \sim V_i} d_{ij}\ge \frac{(c + \mu -2\varepsilon -d )nm}{m^2}\ge(c + \mu -2\varepsilon -d)k. \qedhere
\]
\eproof

Consider a family of $t+1$ clusters in $\mathcal{P}$ which are pairwise $\varepsilon$-regular.
The next simple result finds a large subset in each cluster such that every pair of subsets is super-regular.

\begin{proposition}\label{superregular}
Given a constant $\varepsilon>0$ and integers $m,t$ with $t<\frac{1}{2\varepsilon}$, we let $G$ be an $n$-vertex graph and $V_1,V_2,\ldots, V_{t+1}$ be vertex-disjoint subsets each of size $m$ in $G$ such that every pair $(V_i,V_j)$ is $\varepsilon$-regular with density $d_{ij}:=d(V_i,V_j)$. Then there exists for each $i\in[t+1]$ a subset $V'_i\subseteq V_i$ of size at least $(1-t\varepsilon)m$ such that every pair $(V_i',V_j')$ is $(2\varepsilon,d_{ij}-(t+1)\varepsilon)$-super-regular.
\end{proposition}

\bproof
Fix a pair $(V_i,V_j)$, the definition of regularity implies that there are at most $\varepsilon m$ vertices $v\in V_i$ such that $d_{V_j}(v)\le (d_{ij}-\varepsilon)m$. Let $V'_i$ be a subset of $V_i$ such that for any $j\neq i$, each vertex in $V'_i$ has at least $(d_{ij}-\varepsilon)m$ neighbors in $V_j$. Then we have that $|V'_i|\ge m-t\varepsilon m>\frac{1}{2}m$ for each $i\in [t+1]$. Therefore Fact~\ref{slicing lem} implies that every pair $(V'_i,V'_j)$ is $2\varepsilon$-regular. Since each vertex $v$ in $V'_i$ satisfies
\[d_{V_j'}(v)\ge d_{V_j}(v)-t \varepsilon m\ge(d_{ij}-(t+1)\varepsilon)|V'_j|,\]
we have that $(V_i',V_j')$ is $(2\varepsilon,d_{ij}-(t+1)\varepsilon)$-super-regular.
\eproof

\section{Finding almost perfect tilings}\label{sec3}

Our main tools for embedding $K_{r}$'s are Szemer\'edi's Regularity Lemma~\cite{komlos1996} and a result of Koml\'os~\cite{komlos2000tiling}, of which we only state a special case for $K_{a,b}$ as follows.

\btheorem[Koml\'os~\cite{komlos2000tiling}]\label{Komlos thm}
Let $a,b\in\mathbb{N}$ such that $a\le b$. For any $\gamma>0$, there exists an integer $n_{0}=n_{0}(\gamma,a,b)$ such that every graph $G$ of order $n \ge n_{0}$ with $\delta(G) \ge \frac{a}{a+b}n$ contains a $K_{a,b}$-tiling covering all but at most $\gamma n$ vertices.
\etheorem

\bproof[Proof of Lemma~\ref{tiling lem}]
Given $r,\ell\in \mathbb{N}$ such that $r>\ell\ge\frac{3}{4}r$ and $\mu>0,\delta>0$, we choose
\[
\frac{1}{n} \ll \alpha\ll \frac{1}{N} \ll \varepsilon \ll\mu,\delta,r,\ell.
\]
Let $G$ be an $n$-vertex graph with $\delta(G) \ge \left( \frac{1}{2} + \mu \right)n $ and $\alpha_{\ell}(G) \le \alpha n$. Then by Lemma~\ref{reg} applied on $G$ with the constants $\varepsilon $ and $d:=\frac{\mu}{2}$, we obtain a partition $\mathcal{P}=\{V_{0},\ldots,V_{k}\}$ for some $\frac{1}{\varepsilon}\le k \le N$ and a spanning subgraph $G' \subseteq G$ with properties (a)-(e) as stated.
Let $d':=\frac{r-\ell-1 + \mu}{\ell}$ and $R_{d'}$ be the $d'$-reduced graph of $\mathcal P$.
In particular, $d'<\frac{1}{2}$ as $\ell \ge \frac{3r}{4}$.

\bclaim
$\delta(R_{d'}) \ge \frac{r-\ell}{r} k$.
\eclaim
\bproof
We deduce this by double counting. Fix an arbitrary $V_i$ and simply write $d_R(V_i)$ for the degree of $V_i$ in $R_{d'}$. On one hand, it is easy to observe in $G'$ that
\[
e_{G'}(V_{i},\cup_{j \ne i}V_{j})\ge|V_{i}|(\delta(G')-|V_0|) \ge m \left( \frac{1}{2} + \mu - (d+2\varepsilon ) \right) n \ge \frac{1}{2}k m^{2},
\] where the last inequality follows since $\varepsilon\ll \mu$.
On the other hand, letting $x:=\frac{d_R(V_i)}{k}$, we have that
\[e_{G'}(V_{i},\cup_{j \ne i}V_{j})\le
xkm^{2} + (1-x)kd'm^{2}=x(1-d')km^{2} + d' k m^{2}.
\]
Combining these two bounds, we get that $x\ge\frac{\frac{1}{2} -d'}{1-d'}$. Therefore
\[
x\ge\frac{\frac{1}{2} -d'}{1-d'}=1-\frac{1}{2(1-d')}\ge1-\frac{1}{2\left(1-\frac{r-\ell}{\ell}\right)}\ge1-\frac{\ell}{r},
\]
where the penultimate inequality follows since $d' = \frac{r-\ell-1 + \mu}{\ell}\le\frac{r-\ell}{\ell}$ and in the last inequality we use $1-\frac{r-\ell}{\ell}\ge\frac{r}{2\ell}$ because $\ell \ge \frac{3}{4}r$.
\eproof
Now we apply Theorem~\ref{Komlos thm} on $R_{d'}$ with $\gamma=\frac{\delta}{2},a=r-\ell$ and $b=\ell$ (using that $k\ge \frac{1}{\varepsilon}$ is large enough) to guarantee a family $\mathcal{H}$ of vertex-disjoint copies of $K_{r-\ell,\ell}$ that together cover all but at most $\frac{\delta}{2} k$ vertices in $R_{d'}$.
For each copy of $K_{r-\ell,\ell}$ in $\mathcal{H}$, we shall show that we can find a $K_r$-tiling covering all but at most $2r\ell\varepsilon m$ vertices in (the union of) its clusters.
This would finish the proof as then the union of these $K_r$-tilings would leave at most
\[
|V_0|+\frac{\delta}{2} k m+ |\mathcal{H}|2r\ell\varepsilon m< \varepsilon  n+ \frac{\delta}{2}n+ 2r\ell\varepsilon n< \delta n
\]
vertices uncovered.

Now consider a copy of $K_{r-\ell,\ell}$ in $\mathcal{H}$, and for simplicity we may assume that it has vertex set $\{V_{1},\ldots, V_{r}\}$ such that each pair $(V_i,V_j)$ with $i\in[r-\ell]$ and $j\in[r-\ell+1,r]$ is $\varepsilon$-regular with density at least $d'$.
For each $i\in[r-\ell]$ we partition $V_{i}$ into $\ell$ parts $V_{i,1},\ldots,V_{i,\ell}$ of equal size and each $V_{j}$ with $j\in[r-\ell+1,r]$ into $(r-\ell)$ parts $V_{j,1},\ldots,V_{j,r-\ell}$ of equal size.
Here we may further assume that $|V_{i,i'}|=\frac{m}{\ell}$ for $i\in[r-\ell],i'\in[\ell]$ and $|V_{j,j'}|=\frac{m}{r-\ell}$ for every $j\in[r-\ell+1,r],j'\in[r-\ell]$.
\bclaim\label{embedding lem}
Given $i\in[r-\ell],i'\in[\ell],j\in[r-\ell+1,r]$ and $j'\in[r-\ell]$, $G[V_{i,i'}\cup V_{j,j'}]$ admits a $K_r$-tiling covering all but at most $\frac{2r}{r-\ell}\varepsilon m$ vertices in $V_{i,i'}\cup V_{j,j'}$.
\eclaim
We postpone its proof to the end of the section.
Note that we can arbitrarily pair up all $V_{i,i'}$'s and $V_{j,j'}$'s and then apply Claim~\ref{embedding lem} to each such pair, which together give a $K_r$-tiling covering all but at most $2r\ell\varepsilon m$ vertices in $V_1\cup V_2\cup \cdots\cup V_r$, as desired.
\eproof

\begin{figure}[htb]\label{2}
\center{\includegraphics[width=6.6cm]  {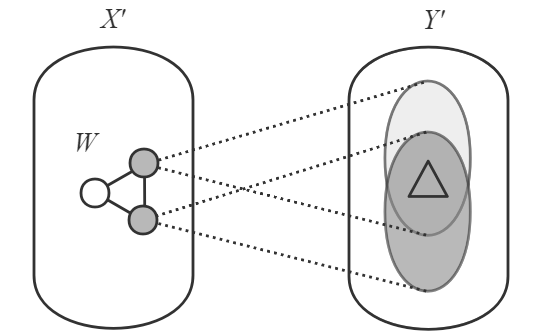}}   \caption{ An illustration of Claim~\ref{embedding lem}.}
\end{figure}
\bproof[Proof of Claim~\ref{embedding lem}]
We shall greedily find vertex-disjoint copies of $K_r$ such that each such copy of $K_r$ contains exactly $r-\ell$ vertices in $V_{i,i'}$ and $\ell$ vertices in $V_{j,j'}$. We write $X:=V_{i,i'}$ and $Y:=V_{j,j'}$ for convenience. Now it suffices to show that for any $X'\subseteq X$ of size at least $2\varepsilon m$ and $Y'\subseteq Y$ of size at least $2\varepsilon m$, there exists a copy of $K_r$ in $X'\cup Y'$ with exactly $r-\ell$ vertices in $X'$ and $\ell$ vertices in $Y'$. Since $(V_i,V_j)$ is $\varepsilon$-regular with density at least $d'$, there exists a subset $X''\subseteq X'$ of size at least $|X'|-\varepsilon m$ vertices such that every vertex in $X''$ has at least $(d'-\varepsilon)|Y'|$ neighbors inside $Y'$. Note that there exists a copy of $K_\ell$ in $X''$ as $\alpha\ll \varepsilon$, whose vertex set is denoted by $W$. Now we have
\[
\sum_{v \in Y'}d_W(v)=e(W,Y')=\sum_{w \in W}d_{Y'}(w) \ge |W|(d'-\varepsilon)|Y'| = |Y'|\left(r-\ell -1 +\frac{\mu}{2}- \varepsilon\ell\right).
\]
By convexity of the function $f(x)=\binom{x}{n}$ where $x\ge n-1$, we have
\[
\sum_{v \in Y'} \binom{d_W(v)}{r-\ell} \ge |Y'| \binom{\frac{1}{|Y'|} \sum_{v \in Y'}d_W(v)}{r- \ell}  \ge 2\varepsilon m\binom{r-\ell -1 +\frac{\mu}{2}- \varepsilon\ell}{r-\ell}.
\]
Since $\frac{\mu}{2}-\varepsilon\ell>\frac{\mu}{4}$, we have
\[
f(\mu):=\binom{r-\ell -1 +\frac{\mu}{2}- \varepsilon\ell}{r-\ell}>\binom{r-\ell-1+\frac{\mu}{4}}{r-\ell}=\frac{(r-\ell-1+\frac{\mu}{4}) \cdots (\frac{\mu}{4})}{(r-\ell)\cdots 1}>\left(\frac{\mu}{4(r-\ell)}\right)^{r-\ell}.
\]
Therefore by the Pigeonhole Principle, there exists a subset $W'$ of $r-\ell$ vertices in $W$ having a set $Y''$ of at least $\frac{2\varepsilon f(\mu) m}{\binom{\ell}{r-\ell}}\ge\frac{2\varepsilon f(\mu)(1-\varepsilon)n}{\binom{\ell}{r-\ell}N}$ common neighbors in $Y'$. Thus we can find a copy of $K_\ell$ in $Y''$ since $\alpha\ll\frac{1}{N}\ll\varepsilon\ll\mu$, which together with $W'$ yields a copy of $K_{r}$ in $G[X' \cup Y']$ as desired.
\eproof

\section{Proof of Lemma~\ref{absorbing lem}}\label{sec4}

\subsection{Lattice-based absorbing method}\label{sub4.1}

To illustrate the lattice-based absorbing method, we need the following notation introduced by Keevash and Mycroft~\cite{keevash2015polynomial}. Let $G$ be an $n$-vertex graph. We will work with a vertex partition $\mathcal{P} = \{V_{1},\ldots,V_{k}\}$ of $V(G)$ for some integer $k \ge 1$. The \emph{index vector} $\mathbf{i}_{\mathcal{P}}(S) \in \mathbb{Z}^{k}$ of a subset $S \subseteq V(G)$ with respect to $\mathcal{P}$ is the vector whose $i$th coordinate is the size of the intersections of $S$ with $V_i$ for each $i\in[k]$. For each $j\in [k]$, let $\mathbf{u}_{j} \in \mathbb{Z}^{k}$ be the $j$th unit vector, i.e. $\mathbf{u}_j$ has 1 on the $j$th coordinate and 0 on the other coordinates. A \emph{transferral} is a vector of the form $\mathbf{u}_{i}-\mathbf{u}_{j}$ for some distinct $i,j \in [k]$. For any $\mathbf{x} = \{x_{1},\ldots,x_{k} \}\in \mathbb{Z}^{k}$, let $|\mathbf{x}| := \sum_{i=1}^{k}x_{i}$. We say that $\mathbf{x} \in \mathbb{Z}^{k}$ is an \emph{$r$-vector} if it has non-negative coordinates and $|\mathbf{x}| = r$. Let $H$ be a graph and let $\beta > 0$. Define $I_{\mathcal{P}}^{\beta}(G)$ to be the set of all $\mathbf{i} \in \mathbb{Z}^{k}$ such that $G$ contains at least $\beta  n$ vertex-disjoint copies of $H$ whose vertex sets all have index vector $\mathbf{i}$ and let $L_{\mathcal{P}}^{\beta }(G)$ denote the lattice (additive subgroup) in $\mathbb{Z}^{k}$ generated by $I_{\mathcal{P}}^{\beta}(G)$.

Another crucial notion of $H$-reachability has been given in the introduction. As mentioned above, we shall construct for every $r$-vertex set a family of linearly many vertex-disjoint absorbers. To achieve this, we make use of the following result from~\cite{han2021ramseyturan}.

\blemma[\cite{han2021ramseyturan}]\label{absorbers}
Given $r,t\in\mathbb{N}$ with $r\ge3$ and $\beta >0$, the following holds for sufficiently large $n\in\mathbb{N}$. Let $G$ be an $n$-vertex graph such that $V(G)$ is $(K_r,\beta,t)$-closed. Then every $S \in \tbinom{V(G)}{r}$ has a family of at least $\frac{\beta}{r^2t} n$ vertex-disjoint $(S,t)$-absorbers.
\elemma

In order to apply Lemma~\ref{absorbers}, we need the following result which guarantees that $V(G)$ is closed.

\blemma\label{merge lem}
Given $r,\ell \in \mathbb{N}$ with $r\ge4$, $r> \ell > \frac{2}{3}r$ and $\mu,\beta_1,\gamma_1$ with $0<\gamma_1,\beta_1<1$, there exist positive constants $\alpha, \beta<\beta_1$ and $t\in \mathbb{N}$ such that the following holds for sufficiently large $n$. Let $G$ be an $n$-vertex graph with $\delta(G) \ge \left( \frac{1}{2} + \mu\right)n $ and $\alpha_{\ell}(G) \le \alpha n$ such that every vertex in $V(G)$ is $(K_r,\beta_1,1)$-reachable to at least $\gamma_1n$ other vertices. Then $V(G)$ is $(K_{r},\beta,t)$-closed.
\elemma

The conditions in Lemma~\ref{merge lem} are too strong to obtain in our graph $G$.
However, we are able to show the following, where a small exceptional set $B$ is allowed and (crucially) the induced subgraph $G-B$ is good enough.

\blemma\label{reachable lem}
Given $r,\ell\in\mathbb{N}$ with $r > \ell \ge 2$ and $\tau, \mu$ with $0<\tau<\mu$, there exist positive constants $\alpha, \beta_1, \gamma_1$ and $n_{0} \in \mathbb{N}$ such that the following holds for every $n\ge n_0$. Let $G$ be an $n$-vertex with $\delta(G) \ge \left( \frac{1}{2-\varrho^*_{\ell}(r-1)} + \mu \right)n $ and $\alpha_{\ell}(G)\le \alpha n$. Then $G$ admits a partition $V(G)=B\cup U$ such that $|B|\le \tau n$ and every vertex in $U$ is $(K_r,\beta_1,1)$-reachable $({\rm in} \ G[U])$ to at least $\gamma_1n$ other vertices in $U$.
\elemma

Note that to invoke Lemma~\ref{merge lem} on $G[U]$, we crucially need that the reachability of $U$ is established in $G[U]$, rather than in $G$.
We also remark that Lemma~\ref{reachable lem} is one of the key places where the minimum degree condition $\delta(G) \ge \left( \frac{1}{2-\varrho^*_{\ell}(r-1)} + \mu \right)n $ is needed.

We shall briefly describe how to use these lemmas. We first use Lemma~\ref{reachable lem} and get a partition $V(G)=B\cup U$ such that every vertex in $U$ is reachable to linearly many other vertices in $U$. Applying Lemma~\ref{merge lem} and Lemma~\ref{absorbers} to $G[U]$, we get the following corollary.

\bcorollary\label{Coro}

Given positive integers $r,\ell$ with $r>\ell>\frac{2}{3}r$ and $\tau, \mu$ with $0<\tau<\mu$, there exist $\alpha>0, \beta < \mu^3$ and $t,n_0\in \mathbb{N}$ such that the following holds for every $n\ge n_0$. Let $G$ be an $n$-vertex graph with $\delta(G) \ge \left( \frac{1}{2-\varrho^*_{\ell}(r-1)} + \mu \right)n $ and $\alpha_{\ell}(G) \le \alpha  n$. Then $G$ admits a partition $V(G)=B\cup U$ such that $|B|\le \tau n$ and every $S \in \tbinom{U}{r}$ has a family of at least $\frac{\beta}{r^2t} n$ vertex-disjoint $(S,t)$-absorbers in $U$.

\ecorollary

Now we are ready to prove Lemma~\ref{absorbing lem} using Corollary~\ref{Coro} and Lemma~\ref{Nenadov lem}.

\bproof[Proof of Lemma~\ref{absorbing lem}]
Given positive integers $\ell,r$ with $r>\ell>\frac{2r}{3}$ and $\mu,\gamma$ with $0<\gamma\le\frac{\mu}{2}$. Choose $\frac{1}{n}<\alpha\ll\xi\ll\beta,\frac{1}{t},\tau\ll \gamma$. Let $G$ be an $n$-vertex graph with $\delta(G) \ge \left( \frac{1}{2-\varrho^*_{\ell}(r-1)} + \mu \right)n $ and $\alpha_{\ell}(G) \le \alpha n$.
Corollary~\ref{Coro} implies that $G$ admits a partition $V(G)=B \cup U$ such that $|B|\le \tau n$ and every $S \in \tbinom{U}{r}$ has a family of at least $\frac{\beta}{r^2t} n$ vertex-disjoint $(S,t)$-absorbers in $U$. Let $G_1:=G[U]$.
Applying Lemma~\ref{Nenadov lem} on $G_1$, we obtain in $G_1$ a $\xi$-absorbing subset $A_1$ of size at most $\frac{\beta}{r^2t} n$.

Now, we shall iteratively pick vertex-disjoint copies of $K_r$ each covering at least a vertex in $B$ whilst avoiding using any vertex in $A_1$.
Let $W:=A_1$ and $G_2:=G-W$.
For $u\in B$, we apply Proposition~\ref{def} iteratively to find a copy of $K_r$ covering $u$ in $G_2$, while avoiding $A_1\cup B$ and all copies of $K_r$ found so far.
This is possible as during the process, the number of vertices that we need to avoid is at most $|A_1|+r|B|\le \frac{\beta}{r^2t}n + r \tau n\le\frac{\mu}{2}n$, by our choice of the constants.
Let $K$ be the union of the vertex sets over all copies of $K_r$ covering $B$ and $A:=A_1\cup K$.
Recall that $A_1$ is a $\xi$-absorbing set for $G_1=G-B$, and $B\subseteq K\subseteq A$.
Thus,  $A$ is a $\xi$-absorbing set for $G$ simply because we do not need to consider the possibility of absorbing vertices in $B$.
Moreover, we have
\[
|A|=|A_1|+|K|\le \frac{\beta}{r^2t} n+r\tau n\le \gamma n,
\]
where the last inequality follows since $\beta,\frac{1}{t},\tau\ll\gamma$.
The proof is completed.
\eproof

The rest of this paper is devoted to the proofs of Lemmas~\ref{merge lem} and~\ref{reachable lem}.

\section{Proof of Lemma~\ref{reachable lem}}\label{sec5}

Recall that given our minimum degree condition $\delta(G) \ge \left( \frac{1}{2-\rho(\ell, r-1)} + \mu \right)n $, Proposition~\ref{def} shows that for every $v$ in $G$, one can embed linearly many vertex-disjoint copies of $K_{r-1}$ in $N(v)$.
What we need in Lemma~\ref{reachable lem} is that after excluding $o(n)$ exceptional vertices, for every vertex $v$ in $G$, there exists at least $\Omega(n)$ vertices $w$ such that one can embed linearly many vertex-disjoint copies of $K_{r-1}$ in $N(v)\cap N(w)$.
This is a considerable strengthening and its proof is certainly much more involved.
We use the regularity method and then naively, when we fix $v\in V_i$ for some $i$ in the regular partition, we tend to let $w$ be in the same cluster.
In fact, for most $w\in V_i$, $|N(w)\cap N(v)\cap V_j|$ should be as expected as long as $d_{ij}=d(V_i, V_j)$ is not too small.
Then we tend to find a suitable induced subgraph on which we could establish a minimum degree condition so that we could use the definition of $\rho(\ell, r-1)$, that is, proceed as in the proof of Proposition~\ref{def}.
There are also some technical issues that we must deal with, for instance, we can not work with all clusters in the reduced graph at the same time as there would be too many exceptional vertices so that we can not establish a minimum degree condition.
In the actual proof we shall partition the reduced graph randomly into blocks of constant size so that we can focus on a constant number of clusters where we could safely ignore vertices with wrong degrees and establish super-regularity.
It is routine to show that each part of the random partition inherits the minimum degree condition, but we actually need to pass down a somewhat stronger degree-estimate (see~\eqref{ineq1}).

We start with the following result, which shows that we can find a family of vertex-disjoint subsets of $R_d^*$ inheriting the degree-estimate in~\eqref{ineq1} in $R_d^*$.
We state this in a slightly more general form, that is, for edge-weighted graphs.

\blemma\label{random partition}
Let $\mu,c$ be any positive constants. Then there exist integers $t_0=t_0(\mu)$ and $k_0=k_0(\mu)$ such that the following holds for $k\ge k_0$ and $t\ge t_0$. Suppose $R$ is a $k$-vertex edge-weighted graph on vertex set $\{w_1,\ldots,w_k\}$ and let $d_{ij}$ denote the weight of $w_iw_j\in E(R)$ such that for any $i,j\in[k]$
\begin{equation}\label{ineq1}
\sum_{p\in[k]\setminus\{i,j\}}\left(d_{ip}d_{jp}-\left(c+\frac{\mu}{6}\right)d_{ip}\right) \ge \frac{\mu}{24}k.
\end{equation}
Then there exists a family of vertex-disjoint $(t+1)$-sets $S_1,\ldots,S_Q$ with $Q\ge\left(1-\frac{\mu}{8(t+1)}\right)\left\lfloor\frac{k}{t+1}\right\rfloor$ such that for each $q\in [Q]$, $R[S_q]$ has the following property: for any $w_i,w_j\in S_q$,
\begin{equation}\label{E}
  \sum_{w_p \in S_q\setminus \{w_i,w_j\}}\left(d_{ip}d_{jp}-\left(c+\frac{\mu}{6}\right)d_{ip}\right)>0.
\end{equation}
\elemma
The proof of Lemma~\ref{random partition} is based on Lemma~\ref{random set} which shows that random subgraphs of weighted graphs almost inherit condition \eqref{ineq1}. We postpone its proof to the end of the section and use it to prove Lemma~\ref{reachable lem} first.

\bproof[Proof of Lemma~\ref{reachable lem}]
Given $r,\ell\in \mathbb{N}$ with $r>\ell\ge2$ and $\tau,\mu>0$, write $\rho:=\varrho^*_{\ell}(r-1)$.
Choose
\[
\frac{1}{n}\ll \alpha\ll\beta_1,\gamma_1\ll \frac{1}{N}\ll \varepsilon\ll\frac{1}{t}\ll\tau, \mu, \rho
\]
and fix an $n$-vertex graph $G$ with $\delta(G) \ge \left( \frac{1}{2-\rho} + \mu \right)n $ and $\alpha_{\ell}(G)\le \alpha n$.
By applying Lemma~\ref{reg} on $G$ with $\varepsilon>0$ and $d:=\frac{\mu}{4}$, we obtain an $\varepsilon$-regular partition $\mathcal{P}=V_{0}\cup \cdots \cup V_{k}$ of $V(G)$ for some $\frac{1}{\varepsilon} \le k \le N$ and a spanning subgraph $G'$ of $G$ with properties (a)-(e) as stated.
Let $R^*_d$ be the weighted reduced graph on $\{ V_{1}, \ldots , V_{k} \}$ as aforementioned, where we use $d_{ij}$ to denote the density of $(V_i, V_j)$ in $G'$.
\bclaim\label{property}
$R^*_d$ satisfies inequality \eqref{ineq1} with $c=\rho$.
\eclaim
We shall prove Claim~\ref{property} after the proof of Lemma~\ref{reachable lem} and here we take Claim~\ref{property} for granted. Then applying Lemma~\ref{random partition} on $R=R^*_d$ with $c=\rho$, we obtain a family $\mathcal{S}$ of vertex-disjoint $(t+1)$-sets $S_1,\ldots,S_Q$ in $R^*_d$ for some $Q\ge\left(1-\frac{\mu}{8(t+1)}\right)\left\lfloor\frac{k}{t+1}\right\rfloor$ such that for each $q\in [Q]$, $R^*_d[S_q]$ satisfies \eqref{E}.
Next we find the partition $V(G)=B\cup U$.

Let $B_1\subseteq V(G)$ be the union of all clusters which are not included in any $S_i$ with $i\in[Q]$. By renaming if necessary, let $V_1, V_2,\ldots,V_{k'}$ be the remaining clusters with $V_1\cup\cdots\cup V_{k'}=V(G)\setminus~B_1$ where $k'\ge\left(1-\frac{\mu}{8(t+1)}\right)k$. Fix $S_q$ with $q\in[Q]$, without loss of generality, we write $S_q=\{V_1,V_2,\ldots,V_{t+1}\}$. Then applying Proposition~\ref{superregular} on $G$ with $\varepsilon\ll \frac{1}{t}$ and letting $m':=(1-t\varepsilon)m$, we obtain for each $i\in[t+1]$ a subset $V'_i\subseteq V_i$ of size at least $m'$ such that every pair $(V_i',V_j')$ is $(2\varepsilon,d_{ij}-(t+1)\varepsilon)$-super-regular.
Here we define $B^q:=\bigcup_{i\in[t+1]}(V_i\setminus V_i')$ and let $B_2$ be the union of all $B^q$'s with $q\in[Q]$.
Finally we claim that $B:=B_1 \cup B_2$ is desired, in particular, $V_0\subseteq B_1$.
Note that $B$ has size
\[
|B_1|+|B_2|\le|V_0| +(k-Q(t+1))m + k'\varepsilon tm \le \varepsilon n + \frac{\mu n}{8(t+1)} + \varepsilon t n <\tau n,
\]
as $\varepsilon\ll\frac{1}{t}\ll\tau$. 
Now it suffices to show that for each $i\in[k']$, every vertex in $V'_i$ is $(K_r,\beta_1,1)$-reachable to many other vertices in $V'_i$. Claim~\ref{claim5.3} below finishes the proof of Lemma~\ref{reachable lem}.
\eproof

It remains to prove Claim~\ref{property} and the following claim.

\bclaim\label{claim5.3}
For any $i\in[k']$, every vertex in $V'_i$ is $(K_r,\beta_1,1)$-reachable to at least $\gamma_1 n$ vertices in $V'_i$.
\eclaim
\bproof
Without loss of generality, suppose $S_1=\{V_1,\ldots,V_{t+1}\}$ and  we prove the claim for $V'_1$.
Let $I_1$ be the set of all indices $i\in[t+1]$ such that $V_i$ is adjacent to $V_1$ in $R^*_d$, that is $d_{1i}>0$.
Note that if $I_1= \emptyset$, then all $d_{1i}=0$, which contradicts~\eqref{E}.
Recall that for any $i\in I_1$, $(V_1',V_i')$ is $(2\varepsilon,d_{1i}-(t+1)\varepsilon)$-super-regular with $d_{1i}\ge d$.
Then every $u\in V'_1$ has at least $(d_{1i} - \varepsilon (t+1) )|V'_{i}|\ge(d-\varepsilon  (t+1))|V'_i|\ge\frac{\mu}{8} m'$ neighbors in $V'_{i}$.
Given a vertex $u\in V'_1$, we choose for each $i\in I_1$ a set $N_i(u)\subseteq N_{G'}(u)\cap V'_i$ of size exactly $m_0$ with $m_0:=\frac{\mu}{8}m'$.

Now we pick a subset $V^{''}_{1}\subseteq V'_1$ such that for each $i\in I_1$, every $v \in V''_1$ satisfies \[(d_{1i}-\varepsilon)m_0\le|N(v)\cap N_{i}(u)|\le(d_{1i}+\varepsilon)m_0.\]
We claim that $u$ is $(K_r,\beta_1,1)$-reachable to every $v \in V_1''$. This essentially completes the proof since \[|V_1''|\ge m'-2t\varepsilon m\ge (1-3t\varepsilon)m\ge \gamma_1 n,\] where the first inequality follows from the fact that $(V_1,V_i)$ is $\varepsilon$-regular for each $i\in I_1$, while the last inequality follows since $\gamma_1\ll\frac{1}{N} \ll \varepsilon$. Now it remains to show that for every $v\in V_1''$, $u$ and $v$ are $(K_r,\beta_1,1)$-reachable.
The key step is to build a subgraph with a large minimum degree condition, as mentioned in the proof outline at the beginning of this section.

\bclaim
There exists a subset $T'$ in $\bigcup_{i\in I_1}(N_{i}(u)\cap N(v))$ such that $\delta(G[T'])\ge\left(\rho+\frac{\mu}{12}\right)|T'|$.
\eclaim
For simplicity, we let $N_{i}:=N_{i}(u)\cap N(v)$ for every $i\in I_1$ and $T:= \bigcup_{i\in I_1}N_{i}$. Note that for every $i\in I_1$, we have
\begin{equation}\label{ineq2}
(d-\varepsilon)\frac{\mu}{9}m \le (d_{1i}-\varepsilon)m_0\le|N_{i}|\le(d_{1i}+\varepsilon)m_0.
\end{equation}
Choose a subset $N^*_i\subseteq N_{i}$ such that for any $j\in I_1$ with $j\neq i$, each vertex in $N^*_i$ has at least $(d_{ij}-\varepsilon)|N_j|$ neighbors in $N_j$. Recall that $|N_i|\ge (d_{1i}-\varepsilon)m_0> \varepsilon m$ since $\varepsilon\ll \mu$. Then using the fact that $(V_i,V_j)$ is $\varepsilon$-regular, we have that $|N^*_{i}|\ge |N_{i}|-t\varepsilon m$ for each $i\in I_1$ and every vertex $v$ in $N^*_{i}$ has \[d_{N^*_j}(v)\ge d_{N_j}(v)-t \varepsilon m\ge(d_{ij}-\varepsilon)|N_j|-t\varepsilon m.\]
Let $T':=\bigcup_{i\in I_1}N^*_i$. Then
\begin{equation}\label{ineq3}
\sum_{i\in I_1}(|N_i|-t\varepsilon m)\le|T'|\le \sum_{i\in I_1}|N_i|\le \sum_{i\in I_1}(d_{1i}+\varepsilon)m_0.
\end{equation}
In particular, as $I_1\neq \emptyset$ we get $|T'|\ge (d-\varepsilon)\frac{\mu}{9}m - t\varepsilon m \ge \frac{d\mu}{10} m \ge\frac{d\mu}{11k} n$.
For every $w\in N^*_j$ with $j\in I_1$, we have
\begin{align}
d_{T'}(w)\ge\sum_{p\in I_1\setminus\{j\}}d_{N^*_p}(w) \nonumber
&\ge\sum_{p\in I_1\setminus\{j\}}(d_{jp}-\varepsilon)|N_p|-t^2\varepsilon m  \\  \nonumber
\underrightarrow{\eqref{ineq2}}~~~~~&\ge \sum_{p\in I_1\setminus\{j\}}(d_{jp}-\varepsilon)(d_{1p}-\varepsilon)m_0-t^2\varepsilon m \\  \nonumber
&\ge \sum_{p\in I_1\setminus\{j\}}d_{jp}d_{1p}m_0-3t^2\varepsilon m   \\  \nonumber
\underrightarrow{\eqref{E}}~~~~~ &\ge \left(\rho+\frac{\mu}{6}\right)\sum_{p\in I_1\setminus\{j\}}d_{1p}m_0-3t^2\varepsilon m \\ \nonumber
\underrightarrow{\eqref{ineq3}}~~~~~&\ge\left(\rho+\frac{\mu}{12}\right)|T'|,  \nonumber
\end{align}
where the last inequality follows from \eqref{ineq3} and $\varepsilon\ll \frac{1}{t}\ll\mu$.

So far we obtain that $G[T']$ is an induced subgraph of $G$ with $\delta(G[T'])\ge\left(\rho+\frac{\mu}{12}\right)|T'|$ and $\alpha_\ell(G[T'])\le\alpha n \le \frac{11\alpha k}{d\mu}|T'|$. Let $\mathcal{K}$ be a maximal family of vertex-disjoint copies of $K_{r-1}$ in $T'$. Then we claim that $|\mathcal{K}|\ge \beta_1 n$.
Suppose for contrary that $|\mathcal{K}|< \beta_1 n$ and let $K$ be the union of the vertex sets over all copies of $K_{r-1}$ in $\mathcal{K}$. Then since $\beta_1\ll\varepsilon$, we obtain that $|K|<r\beta_1 n< \frac{\mu}{24}|T'|$, and thus
$\delta(G[T'\setminus K])\ge\left(\rho+\frac{\mu}{24}\right)|T'|$.
By the definition of $\varrho^*_{\ell}(r-1)$ and the fact that $\alpha \ll \frac{1}{N} \ll\mu$, we can pick a copy of $K_{r-1}$ in $T'\setminus K$, a contradiction.

Thus by definition, we infer that $u$ and $v$ are $(K_r,\beta_1,1)$-reachable, which completes the proof.
\eproof
Now we prove Claim~\ref{property}.
\bproof[Proof of Claim~\ref{property}]
Recall that we have an $\varepsilon$-regular partition $\mathcal{P}=V_{0}\cup \cdots \cup V_{k}$ of $V(G)$ for some $\frac{1}{\varepsilon}\le k \le N$ and a spanning subgraph $G'$ of $G$ with properties (a)-(e) as stated.
Consider the graph $G''=G'- V_0$ and fix two clusters $V_i$ and $V_j$. By Fact~\ref{min degree}, we have $\delta(R^*_d)\ge\left( \frac{1}{2-\rho} + \frac{\mu}{2} \right)k$. Now we want to count in $G''$ the number $S_{ij}$ of all copies of $P_3$ which start from $V_i$, pass through $V_p$ and end in $V_j$ for all $p\in[k]\setminus\{i,j\}$. For the upper bound, since $(V_i,V_p)$ is $\varepsilon $-regular, it is easy to see that at least $(m-\varepsilon  m)$ vertices in $V_i$ have degree at most $(d_{ip}+\varepsilon )m$ in each $V_p$. Then
\[
\begin{split}
S_{ij}&\le\sum_{p \in [k]\setminus\{i,j\}}\left(|V_i|(d_{ip}+\varepsilon )m(d_{jp}+\varepsilon )|V_j|+\varepsilon m \cdot e(V_j,V_p)\right)\\
&\le m^3 \sum_{p \in [k]\setminus\{i,j\}}(d_{ip}+\varepsilon )(d_{jp}+\varepsilon )+ \varepsilon  m^3 \sum_{p \in [k]\setminus\{i,j\}}d_{jp}\le m^3 \sum_{p \in [k]\setminus\{i,j\}}d_{ip}d_{jp}+4m^3\varepsilon k,
\end{split}
\]
where in the last inequality we use $\sum_{p\in[k]\setminus\{i,j\}}d_{ip}< k$. For the lower bound, we have that
\[
\begin{split}
S_{ij}&\ge\sum_{v\in V_j}\sum_{u\in V_i}(d_{G''}(u)+d_{G''}(v)-n-2m) \ge \sum_{v\in V_j}\sum_{u\in V_i}(d_{G''}(u)+\delta(G')-|V_0|-n-2m)\\
&\ge m^3\left(\sum_{p \in [k]\setminus\{i,j\}}d_{ip}+\left(\frac{1}{2-\rho} + \frac{\mu}{2}\right)k-k\right)\\
&=m^3\left(\left(\rho+\frac{\mu}{6}+1-\rho-\frac{\mu}{6}\right)\sum_{p \in [k]\setminus\{i,j\}}d_{ip}+\left(\frac{1}{2-\rho} + \frac{\mu}{2}\right)k-k\right)\\
&\ge m^3\left(\rho+\frac{\mu}{6}\right)\sum_{p \in [k]\setminus\{i,j\}}d_{ip}+\frac{\mu}{12}m^3k,
\end{split}
\]
where in the last inequality we use $\sum_{p\in[k]\setminus\{i,j\}}d_{ip}\ge\delta(R^*_d)\ge\left(\frac{1}{2-\rho} + \frac{\mu}{2} \right)k$.
Combining the upper bound and the lower bound, we have
\[
\sum_{p \in [k]\setminus\{i,j\}}d_{ip}d_{jp}+4\varepsilon k\ge \left(\rho+\frac{\mu}{6}\right)\sum_{p \in [k]\setminus\{i,j\}}d_{ip}+\frac{\mu}{12}k.
\]
Therefore by the fact that $\varepsilon\ll \mu$, we get
\[
\sum_{p \in [k]\setminus\{i,j\}}d_{ip}d_{jp}-\left(\rho+\frac{\mu}{6}\right)\sum_{p \in [k]\setminus\{i,j\}}d_{ip}\ge\frac{\mu}{24}k.
\]
This completes the proof.
\eproof

To prove Lemma~\ref{random partition}, we need the following concentration result~\cite[Corollary 2.2]{greenhill2017average}.

\blemma[\cite{greenhill2017average}]\label{con ineq}
Let $\tbinom{[N]}{r}$ be the set of $r$-subsets of $\{ 1,\ldots,N\}$ and let $h:\tbinom{[N]}{r} \to \mathbb{R}$ be given. Let $C$ be a uniformly random element of $\tbinom{[N]}{r}$. Suppose that there exists $c' \ge 0$ such that
\[
   |h(A) - h(A')| \le c'
\]
for any $A,A'\in \tbinom{[N]}{r}$ with $|A \cap A'| = r -1 $. Then
\[
   \mathbb{E}e^{h(C)} =\exp(\mathbb{E}h(C) + K)
\]
where $K$ is a real constant such that $0\le K \le \frac{1}{8}{\rm min}\{r,N-r\} c'^{2}$. Furthermore, for any real $m > 0$,
\[
  \mathbb{P}[|h(C)-\mathbb{E}h(C)| \ge m] \le 2 \exp \left(- \frac{2m^{2}}{{\rm min}\{r,N-r\} c'^{2}}\right).
\]
\elemma

We also make use of the following result whose proof is based on Lemma~\ref{con ineq}.
\blemma\label{random set}
Given positive constants $c,\mu$ and $k,t\in \mathbb{N}$, we fix a $k$-vertex edge-weighted graph $R$ on vertex set $W=\{w_1,\ldots,w_k\}$ satisfying inequality~\eqref{ineq1}. Let $S$ be a $(t+1)$-set chosen from $\binom {W}{t+1}$ uniformly at random. Then with probability at least $1-\binom{t+1}{2}\exp\left(-\frac{\mu^{2}}{1000}t\right)$, the following event happens: for any $w_i,w_j\in S$, it holds that
\begin{align}
\sum_{w_p \in S\setminus\{w_i,w_j\}}\left(d_{ip}d_{jp}-\left(c+\frac{\mu}{6}\right)d_{ip}\right) >0.
\tag{\ref{E}A}
\end{align}
\elemma

\bproof
For simplicity, we use $E$ to denote the event given as above.
For all $i,j\in[k]$, we write $A_{i,j}$ for the event that $\{w_i,w_j\}\subseteq S$, and $B_{i,j}$ for the event that
\begin{align}
\sum_{w_p \in S\setminus\{w_i,w_j\}}\left(d_{ip}d_{jp}-\left(c+\frac{\mu}{6}\right)d_{ip}\right) \le0. \nonumber
\end{align}
Then it holds that $\overline{E}=\bigcup\limits_{1\le i<j\le k}(A_{i,j}\cap B_{i,j})$. Since $\mathbb{P}[A_{i,j}]=\frac{\binom{k-2}{t-1}}{\binom{k}{t+1}}=\frac{\binom{t+1}{2}}{\binom{k}{2}}$, we have that
\begin{align}
\mathbb{P}\left[\bigcup_{1\le i<j\le k}(A_{i,j}\cap B_{i,j})\right]
\le \sum_{1\le i<j\le k}\mathbb{P}[A_{ij}\cap B_{ij}]
= \frac{\binom{t+1}{2}}{\binom{k}{2}}\sum_{1\le i<j\le k}\mathbb{P}[B_{i,j}|A_{i,j}].\nonumber
\end{align}
Now it suffices to prove for all $i,j\in[k]$ that $\mathbb{P}[B_{i,j}|A_{i,j}]\le \exp(-\frac{\mu^{2}}{1000}t)$. Here by symmetry we may show that it holds for $i=1$ and $j=2$. Fix two vertices $w_{1},w_{2}$ in $S$, and order the remaining vertices of $R$ as $w_{3},\ldots,w_{k}$ uniformly at random. Now, condition on $w_1,w_2\in S$, we take $S' = \{w_{3}, \ldots, w_{t+1}\}$. Note that $S'$ is a subset of size $t-1$ chosen from $R\setminus\{w_1,w_2\}$ uniformly at random. Let $f(S'):=\sum_{w_p\in S'}\left(d_{1p}d_{2p}-(c+\frac{\mu}{6})d_{1p}\right)$. Then $\mathbb{P}[B_{1,2}|A_{1,2}]=\mathbb{P} [ f(S')\le 0 ]$.

Let $X_{p}$ be the indicator variable for the event $\{w_p \in S'\}$ for every $p\in[3,k]$. Then $\mathbb{E}[X_p]=\mathbb{P}[w_p \in S']=\frac{t-1}{k-2}.$ Thus we have $f(S')=\sum_{w_p\in [3,k]}X_p\left(d_{1p}d_{2p}-(c+\frac{\mu}{6})d_{1p}\right)$ and
\[
\begin{split}
\mathbb{E}\big{[}f(S')\big{]}&=\mathbb{E}\bigg{[}\sum_{p\in [3,k]}X_p\left(d_{1p}d_{2p}-\left(c+\frac{\mu}{6}\right)d_{1p}\right)\bigg{]}\\
&=\sum_{p\in [3,k]}\left[ \left(d_{1p}d_{2p}-\left(c+\frac{\mu}{6}\right)d_{1p} \right) \mathbb{E}[X_p]\right]\\
\underrightarrow{\eqref{ineq1}}~~&\ge \frac{\mu}{24}k\cdot\frac{t-1}{k-2}> \frac{\mu}{24}(t-1).
\end{split}
\]
A crucial observation is that swapping any $w_j\in S'$ with any $w_i\in V(R)\setminus S$ affects $f(S')$ by at most $1$, which is because $-(c+\frac{\mu}{6})\le d_{1p}d_{2p}-(c+\frac{\mu}{6})d_{1p}\le1-(c+\frac{\mu}{6})$. By Lemma~\ref{con ineq} with $c'=1,m=\mathbb{E}\left[f(S')\right]$, we know that
\[
\mathbb{P}\left[f(S') \le 0 \right] \le 2 {\rm exp} \left(- \frac{2(\frac{\mu}{24}(t-1))^2}{(t-1) \cdot 1^2} \right) < \exp \left(-\frac{\mu^{2}}{1000}t \right).
\]
This completes the proof.
\eproof

Now, we are ready to prove Lemma~\ref{random partition}.
\bproof[Proof of Lemma~\ref{random partition}]
Choose $t_0 \in \mathbb{N}$ such that $8t^3<\mu e^{\frac{\mu^{2}}{1000}t}$ holds for every $t\ge t_0$. Given integers $t\ge t_0$ and $Z=\left\lfloor\frac{k}{t+1}\right\rfloor$, we uniformly and randomly choose from $V(R)$ vertex-disjoint sets $S_1,\ldots,S_Z$ with $|S_q|=t+1$ for each $q\in [Z]$. We say $S_q$ is \emph{good} if the induced subgraph $R[S_q]$ satisfies \eqref{E}, otherwise it is \emph{bad}. By Lemma~\ref{random set}, $\mathbb{P}[S_q \ {\rm is} \ {\rm bad}]\le \binom{t+1}{2}\exp\left(-\frac{\mu^{2}}{1000}t\right)$ for each $q\in[Z]$. Let $X$ be the random variable counting the number of bad sets. Then we have
\[
\mathbb{E}[X]=\sum_{q\in[Z]}\mathbb{P}[S_q \ {\rm is} \ {\rm bad}]\le Z\binom{t+1}{2}\exp\left(-\frac{\mu^{2}}{1000}t\right).
\]
By Markov's inequality, we obtain that
\[
\mathbb{P}\left[X\ge\frac{\mu Z}{8(t+1)}\right]\le\frac{8(t+1)}{\mu Z}\mathbb{E}[X]\le\frac{8(t+1)}{\mu}\binom{t+1}{2}\exp\left(-\frac{\mu^{2}}{1000}t\right)<1.
\]
Thus there exists a family of vertex-disjoint good $(t+1)$-sets $S_1,\ldots,S_Q$ with $Q\ge(1-\frac{\mu}{8(t+1)})\left\lfloor\frac{k}{t+1}\right\rfloor$.
\eproof

\section{Proof of Lemma~\ref{merge lem}}\label{sec6}

Here is a brief proof outline.
In order to show that $V(G)$ is closed, a known way is to first partition $V(G)$ into a constant number of blocks each of which is closed, and then try to merge them into larger (still closed) blocks by analyzing the graph structures.
The following result can be used to construct such a partition (see~\cite[Lemma 4.1]{han2021ramseyturan}).
The key of our proof of Lemma~\ref{merge lem} is in the merging process, where we use the regularity method and an embedding result similar as that of Claim~\ref{embedding lem}.
In particular, such an embedding lemma (Lemma~\ref{KR}) allows us to construct vertex-disjoint copies of $K_r$ with different vertex distributions, which can be used to show the existence of a transferral.

\blemma[Partition lemma~\cite{han2021ramseyturan}]\label{partition lem}
For any positive constants $\gamma_1,\beta_1$ and integer $r \ge 2$, there exist $\beta_2=\beta_2(\gamma_1,\beta_{1},r) > 0$ and $t_2\in\mathbb{N}$ such that the following holds for sufficiently large $n$. Let $G$ be an $n$-vertex graph such that every vertex in $V(G)$ is $(K_{r},\beta_{1},1)$-reachable to at least $\gamma_1 n$ other vertices. Then there is a partition $\mathcal{P} = \{V_{1},\ldots,V_{p}\}$ of $V(G)$ with $p \le \lceil\frac{1}{\gamma_1}\rceil$ such that for each $i \in [p]$, $V_{i}$ is $(K_{r},\beta_2,t_2)$-closed and has size at least $\frac{\gamma_1}{2} n$.
\elemma

The following result~\cite[Lemma 4.4]{han2021ramseyturan} gives a sufficient condition that allows us to iteratively merge two distinct parts into a closed one, given the existence of a transferral.
\blemma[\cite{han2021ramseyturan}]\label{transferral lem} Given any positive integers $r,t$ with $r\ge 3$ and a constant $\beta >0$, the following holds for sufficiently large $n$. Let $G$ be an $n$-vertex graph with a partition $\mathcal{P}=\{V_{1},\ldots,V_{p}\}$ of $V(G)$ such that each $V_{i}$ is $(K_{r},\beta,t)$-closed. For distinct $i,j \in [p]$, if there exist two $r$-vectors $\mathbf{s}, \mathbf{t} \in I_{\mathcal{P}}^{\beta }(G)$ such that $\mathbf{s}-\mathbf{t} = \mathbf{u}_{i} - \mathbf{u}_{j}$, then $V_{i} \cup V_{j}$ is $(K_{r},\frac{\beta}{2} ,2rt)$-closed.
\elemma

To apply Lemma~\ref{transferral lem}, we shall make use of the following lemma to guarantee the existence of a transferral, whose proof is very similar to that of Claim~\ref{embedding lem}.

\blemma\label{KR}
Given $r,\ell\in\mathbb{N}$ with $r\ge\ell\ge\frac{r}{2}\ge2$ and constants $0<4\varepsilon<\mu<1,0<\delta<1$, there exist $\alpha,\gamma$ with $0<\alpha<\gamma$ such that the following holds for sufficiently large $n$. Let $G$ be an $n$-vertex graph with $\alpha_{\ell}(G)\le \alpha n$ and $(V_i,V_j)$ be an $\varepsilon$-regular pair in $G$ with $|V_i|=|V_j|\ge\delta n$. If $d(V_i,V_j) \ge \frac{r- \ell}{\ell} + \mu$, then there exists for each $x\in\{0,1\}$ a family $\mathcal{K}_x$ of at least $\gamma n$ vertex-disjoint copies of $K_r$ such that every $K\in\mathcal{K}_x$ has $|V(K)\cap V_i|=r-\ell+x$ and $|V(K)\cap V_j|=\ell-x$.
\elemma
\bproof
Choose $\frac{1}{n}<\alpha\ll\gamma\ll\varepsilon,\delta$ and let $m:=|V_i|=|V_j|$.
For $x=0,1$, we shall greedily find vertex-disjoint copies of $K_r$ as desired in $G[V_i\cup V_j]$.
Clearly, it suffices to show that for any $X\subseteq V_i, Y\subseteq V_j$ each of size at least $2\varepsilon m$, there exists for each $x\in \{0,1\}$ a copy of $K_r$ inside $X\cup Y$ with exactly $r-\ell+x$ vertices in $X$ and $\ell-x$ vertices in $Y$.
Since $(V_i,V_j)$ is $\varepsilon$-regular with density at least $\frac{r- \ell}{\ell} + \mu$, there exists a subset $X'\subseteq X$ of size at least $|X|-\varepsilon m\ge \varepsilon m$ vertices such that every vertex in $X'$ has at least $(\frac{r- \ell}{\ell} + \mu-\varepsilon)|Y|$ neighbors inside $Y$. Note that by the fact that $\alpha_{\ell}(G)\le \alpha n$ with $\alpha\ll \varepsilon,\delta$, there exists a copy of $K_\ell$ in $X'$, whose vertex set is denoted by $W$. Now we have
\[
\sum_{v \in Y}d_W(v)=e(W,Y)=\sum_{w \in W}d_{Y}(w) \ge |W|\left(\frac{r- \ell}{\ell} + \mu-\varepsilon\right)|Y| = |Y|\left(r-\ell +\mu\ell- \varepsilon\ell\right).
\]
By convexity of the function $f(x)=\binom{x}{n}$ where $x\ge n-1$, we have for each $x\in\{0,1\}$
\[
\sum_{v \in Y} \binom{d_W(v)}{r-\ell+x} \ge |Y| \binom{\frac{1}{|Y|} \sum_{v \in Y}d_W(v)}{r- \ell+x}  \ge 2\varepsilon m\binom{r-\ell +\mu\ell- \varepsilon\ell}{r-\ell+x}.
\]
Since $\mu-\varepsilon>\frac{\mu}{2}$, we have
\[
f(\mu):=\binom{r-\ell +\mu\ell- \varepsilon\ell}{r-\ell+x}>\binom{r-\ell+\frac{\mu\ell}{2}}{r-\ell+x}>\frac{(r-\ell+\frac{\mu\ell}{2}) \cdots (\frac{\mu\ell}{2})}{(r-\ell+1)\cdots 1}>\left(\frac{\mu\ell}{2(r-\ell+1)}\right)^{r-\ell+1}.
\]
Therefore by the Pigeonhole Principle, there exists a subset $W'$ of $r-\ell+x$ vertices in $W$ having a set $Y'\subseteq Y$ of at least $\frac{2\varepsilon f(\mu) m}{\binom{\ell}{r-\ell}}$ common neighbors. Thus we can find a copy of $K_{\ell-x}$ in $Y'$ since $\alpha_{\ell}(G)\le \alpha n$ and $\alpha\ll\varepsilon,\delta$, which together with $W'$ yields a copy of $K_{r}$ in $G[X' \cup Y']$ as desired.
\eproof

Now we are ready to prove Lemma~\ref{merge lem}.
\bproof[Proof of Lemma~\ref{merge lem}]
Given $r,\ell\in\mathbb{N}$ with $r>\ell>\frac{2}{3}r$ and a constant $\mu>0$, we choose \[\frac{1}{n}<\alpha\ll\beta\ll\gamma\ll\varepsilon\ll\beta_2,\frac{1}{t_2}\ll\beta_1,\gamma_1,\mu.\] Fix an $n$-vertex graph $G=(V,E)$ with $\delta(G) \ge \left( \frac{1}{2} + \mu\right)n $ and $\alpha_{\ell}(G) \le \alpha n$ such that every vertex in $V(G)$ is $(K_r,\beta_1,1)$-reachable to at least $\gamma_1n$ other vertices. Then applying Lemma~\ref{partition lem} on $G$ with $\beta_{1},\gamma_1$, we obtain a partition $\mathcal{P}_{0} = \{V_{1},\ldots,V_{p}\}$ for some integer $p \le \lceil\frac{1}{\gamma_1}\rceil$, where each $V_{i}$ is $(K_{r}, \beta_2,t_2)$-closed and $|V_{i}|\ge \frac{\gamma_1 n}{2}$.

We shall proceed by iteratively merging as many distinct parts as possible. As an intermediate step, we choose positive constants $\beta=\beta_{p+1}\ll\beta_p\ll\cdots\ll\beta_3\ll\beta_2$. Our merging procedure goes as follows: at step $s$ ($s\ge1$), anchoring at the current partition $\mathcal{P}=\mathcal{P}_{s-1} = \{V_1,\ldots,V_{p'}\}$ with $p'=p-s+1\in\mathbb{N}$ and each $V_i$ being $(K_r,\beta_{s+1},(2r)^{s-1} t_2)$-closed, if there exist distinct $i,j\in [p']$ and two $k$-vectors $\textbf{s},\textbf{t}\in I^{2\beta_{s+2}}_{\mathcal{P}}(G)$ such that $\textbf{s}-\textbf{t}=\textbf{u}_i-\textbf{u}_j$, then we merge $V_i,V_j$ as a new part. By applying Lemma~\ref{transferral lem} with $\beta=2\beta_{s+2}$ and $t=(2r)^{s-1} t_2$, we obtain that $V_i\cup V_j$ is $(K_r,\beta_{s+2},(2r)^{s} t_2)$-closed. By renaming if necessary, we end the $s$th step with a new partition $\mathcal{P}= \{V_1,\ldots,V_{p'-1}\}$ where each $V_i$ is $(K_r,\beta_{s+2},(2r)^{s} t_2)$-closed. In this way, we continue until the procedure terminates after at most $p-1$ steps.

Suppose we end up with a final partition $\mathcal{P}= \{V_1,\ldots,V_{p'}\}$ for some integer $p'=p-s$ and $0\le s<p$, where each $V_i$ is $(K_r,\beta_{s+2},(2r)^{s} t_2)$-closed. Now we claim that $p'=1$, and this ends the proof by taking $\beta=\beta_{p+1}$ and $t=(2r)^{p-1} t_2$. Assume for a contradiction that $p'>1$, which means the non-existence of $r$-vectors $\textbf{s},\textbf{t}\in I^{2\beta_{s+3}}_{\mathcal{P}}(G)$ such that $\textbf{s}-\textbf{t}=\textbf{u}_i-\textbf{u}_j$ for some distinct $i,j\in [p']$. Then applying Lemma~\ref{reg} on $G$ with constants $\varepsilon \ll \mu$ and $d:=\frac{\mu}{3}$, we obtain a refinement $\mathcal{P}'=\{V_0\}\cup\{V_{i,j}\subseteq V_i:i\in[p'],j\in[k_i]\}$ of $\mathcal{P}$ and a spanning subgraph $G'\subseteq G$ with the properties (a)-(e), where we let $m:=|V_{i,j}|$ for all $i\in[p'],j\in[k_i]$ and $k:=\sum_{i=1}^{p'} k_i$. Let $R^*_d$ be the weighted reduced graph defined on vertex set $\{ V_{i,j}:i\in[p'],j\in[k_i]\}$. Then Fact~\ref{min degree} implies that
\begin{equation}\label{ineq4}
 \delta(R^*_d) \ge  \left( \frac{1}{2} + \frac{\mu}{2} \right)k.
\end{equation}
Recall that $p'>1$. Now, we have the following two claims.
 For any $i,j\in[p'], i'\in[k_i]$ and $j'\in[k_j]$, an edge joining $V_{i,i'}$ and $V_{j,j'}$ is called a \emph{heavy} edge if $d(V_{i,i'}, V_{j,j'}) \ge \frac{1}{2} + \frac{\mu}{2}$.
\begin{claim}\label{cl4.3}
For all distinct $i,j\in [p']$, there is no heavy edge connecting $V_{i,i'}$ and $V_{j,j'}$ for some $i'\in[K_i],j'\in[k_j]$.
\end{claim}
\bproof
We may take $i=1,j=2$ for instance and assume for a contradiction that $d(V_{1,i'}, V_{2,j'}) \ge \frac{1}{2}  + \frac{\mu}{2}$ for some $i'\in[k_1],j'\in[k_2]$. Since $\frac{r-\ell}{\ell}<\frac{1}{2}$, Lemma~\ref{KR} ensures for each $x\in\{0,1\}$ a family $\mathcal{K}_x$ of at least $\gamma n$ vertex-disjoint copies of $K_r$ such that every $K\in\mathcal{K}_x$ has $|V(K)\cap V_{1,i'}|=r-\ell+x$ and $|V(K)\cap V_{2,j'}|=\ell-x$. This gives two $r$-vectors $\textbf{s}=(r-\ell,\ell,0,\ldots,0),\textbf{t}=(r-\ell+1,\ell-1,0,\ldots,0)\in I_{\mathcal{P}}^{2\beta_{s+3}}(G)$ by choosing $\beta_{s+3}\ll \gamma$, which yields a contradiction since $\textbf{s}-\textbf{t}=\textbf{u}_2-\textbf{u}_1$.
\eproof

\begin{claim}\label{cl4.4}
For every $i\in[p']$, every two vertices $V_{i,j},V_{i,j'}$ connected by a heavy edge do not have any common neighbor outside $V_i$ in $R^*_d$.
\end{claim}
\bproof
Without loss of generality, we may take $i=1,j=1,j'=2$ for instance and assume for the contrary that $V_{2,1}$ is adjacent to both $V_{1,1}$ and $V_{1,2}$, where $d(V_{1,1},V_{1,2})\ge\frac{1}{2} + \frac{\mu}{2}$. Since $V_{1,1},V_{1,2}$ and $V_{2,1}$ are pairwise $\varepsilon$-regular, we obtain a set $V_{2,1}'\subseteq V_{2,1}$ of size at least $(1-2\varepsilon)m$ such that each $v\in V_{2,1}'$ has $d_{V_{1,1}}(v),d_{V_{1,2}}(v)\ge (d-\varepsilon)m\ge \frac{\mu m}{4}\ge\frac{\mu(1-\varepsilon)}{4k}n$. Therefore for every vertex $v \in V_{2,1}'$, letting $N_1=N(v)\cap V_{1,1}$ and $N_2=N(v)\cap V_{1,2}$, Lemma~\ref{KR} (applied on the two parts $N_1,N_2$ with $\delta=\frac{\mu(1-\varepsilon)}{4k}$) ensures the existence of $\gamma n$ vertex-disjoint copies of $K_r$ in $N_1\cup N_2$, yielding that $\textbf{s}':=(r,0,\ldots,0)\in I_{\mathcal{P}}^{2\beta_{s+3}}(G)$ by choosing $\beta_{s+3}\ll \gamma$.
Furthermore, by first taking $v\in V'_{2,1}$ and $K_{r-1}$ in $N_1\cup N_2$ by Lemma~\ref{KR}, one can also greedily find $2\beta_{s+3} n$ vertex-disjoint copies of $K_r$ each having index vector $\textbf{t}':=(r-1,1,0,\ldots,0)$ (with respect to $\mathcal{P}$). This yields a contradiction since $\textbf{s}',\textbf{t}'\in I_{\mathcal{P}}^{2\beta_{s+3}}(G)$ and $\textbf{s}'-\textbf{t}'=\textbf{u}_1-\textbf{u}_2$.
\eproof

Now we show how to derive a contradiction and we may assume that $k_1\le \frac{k}{2}$. For $V_{1,1}$, we claim that there exists $V_{j,j'}$ for some $j\in[2,p']$ and $j'\in[k_j]$ such that $V_{1,1}V_{j,j'}\in E(R^*_d)$.
Indeed, otherwise, we have that $ d_{R^*}(V_{1,1})\le 1\cdot k_1\le\frac{k}{2}$, contradicting \eqref{ineq4}. Without loss of generality, we may assume that $V_{1,1}V_{2,1}\in E(R^*_d)$. Let $B_i$, $i\in\{1,2\}$, be the set of vertices in $R^*_d$ which are connected with $V_{i,1}$ via heavy edges. Then by Claim~\ref{cl4.3}, Claim~\ref{cl4.4} and \eqref{ineq4}, we have
\[
   \left(\frac{1}{2} + \frac{\mu}{2}\right)k \le d_{R^*}(V_{1,1}) < 1\cdot|B_1| + \left(\frac{1}{2}+\frac{\mu}{2}\right)(k-|B_1|-|B_2|)
\]
and
\[
\left(\frac{1}{2} + \frac{\mu}{2}\right)k \le d_{R^*}(V_{2,1}) < 1\cdot|B_2| + \left(\frac{1}{2}+\frac{\mu}{2}\right)(k-|B_1|-|B_2|).
\]
Combining these two inequalities, we get
\[
(1+\mu) k< (1+\mu)k-\mu(|B_1|+|B_2|),
\]
a contradiction. The proof is completed.
\eproof

\section{Proof of Proposition~\ref{RTprop}}\label{sec8}

In this section, we shall prove $\varrho^*_\ell(r)=\varrho_\ell(r)$. We firstly need the following results.

\blemma\label{RT1}
Given $r,\ell\in\mathbb{N}$ and positive constants $\varepsilon,\alpha$, if $0<\frac{1}{n}<\frac{1}{n'}\ll\varepsilon$, then we have $\frac{1}{n'}\textbf{RT}^*_{\ell}(n',K_{r},\alpha n')-\varepsilon\le\frac{1}{n}\textbf{RT}^*_{\ell}(n,K_{r},2\alpha n)$.
\elemma

\bproof
Let $\delta=\frac{1}{n'}\textbf{RT}^*_{\ell}(n',K_{r},\alpha n')$. Consider an $n'$-vertex graph $G'$ with $\delta(G')=\delta n'$. We shall define a graph $G$ from $G'$ where every vertex of $G'$ is replaced by an independent set of size either $\left\lceil\frac{n}{n'}\right\rceil$ or $\left\lfloor\frac{n}{n'}\right\rfloor$ proportionally at random so that the resulting graph $G$ has $n$ vertices. Let $p:=\frac{n}{n'}-\left\lfloor\frac{n}{n'}\right\rfloor$. We randomly partition $V(G')$ into two parts $V_1$ and $V_2$ of size $pn'$ and $(1-p)n'$. For every $v\in G'$, we denote by $X_v$ the number of neighbors of $v$ lying in $V_1$. The random variable $X_v$ is hypergeometrically distributed with $\mathbb{E}[X_v]\ge\delta n'\frac{|V_1|}{|V(G')|}=p \delta n'$.
Therefore, by the Chernoff's inequality for the hypergeometric distribution~\cite{Randomgraph2000}, we get
\[
\mathbb{P}[|X_v-\mathbb{E}[X_v]|\ge\varepsilon \mathbb{E}[X_v]]\le2\exp\left(-\frac{\varepsilon^2}{3}\mathbb{E}[X_v]\right)\le2\exp\left(-\frac{\varepsilon^2\delta pn}{3}\right).
\]
Then by the union bound, there exists an outcome such that for every $v\in G'$, we have $ X_v>(1-\varepsilon)p\delta n'$. Let $G$ be the graph obtained from $G'$ by replacing every vertex in $V_1$ $(V_2)$ by an independent set $W_v$ of size $\left\lceil\frac{n}{n'}\right\rceil$ $(\left\lfloor\frac{n}{n'}\right\rfloor)$, and connecting $W_u$ and $W_v$ by a complete bipartite graph if $uv\in E(G')$. Then for every $w\in G$, we have
\begin{align}
d_{G}(w)&=X_v\left\lceil\frac{n}{n'}\right\rceil+(\delta n'-X_v)\left\lfloor\frac{n}{n'}\right\rfloor= X_v\left(\frac{n}{n'}+(1-p)\right)+(d_{G'}(v)-X_v)\left(\frac{n}{n'}-p\right)\nonumber\\
&=d_{G'}(v)\left(\frac{n}{n'}-p\right)+X_v
\ge\delta n'\left(\frac{n}{n'}-p\right)+(1-\varepsilon)p\delta n' =\delta n-\varepsilon p\delta n'\ge(\delta -\varepsilon)n.\nonumber
\end{align}
Then $G$ has minimum degree $(\delta -\varepsilon)n$ and $\alpha_{\ell}(G)\le 2\alpha n$.
\eproof

\blemma\label{RT2}
Let $\frac{1}{n}\ll\eta,\delta<1$ and $\eta<\frac{\delta}{2}$. Suppose that $G$ is an $n$-vertex graph with at least $\delta\binom{n}{2}$ edges. Then there exists $n'\in\mathbb{N}$ with $n'\ge \frac{\eta n}{4}$ such that $G$ has an $n'$-vertex subgraph $G'$ with $\delta(G')\ge(\delta-\eta)n'$.
\elemma
\bproof
Let $G_0:=G$. For $i\ge0$, we recursively define a sequence of graphs $\{G_i\}$ as follows. If $G_i$ has a vertex $v$ of degree less than $(\delta-\eta)|G_i|$, then we let $G_{i+1}:=G_i-\{v\}$. Repeat this process. If we reach $i=n-\eta\frac{n}{4}$, then we get
\begin{align}
e(G_i)&\ge\delta\binom{n}{2}-(\delta-\eta)\sum_{j=1}^i(n-j)
=\delta\binom{n}{2}-(\delta-\eta)\left(ni-\frac{i(i+1)}{2}\right)\nonumber\\
&=\delta\binom{n}{2}-(\delta-\eta)\left((n-i-1)i+\frac{i(i+1)}{2}\right)
\ge\delta\binom{n}{2}-(\delta-\eta)\left(\frac{\eta n}{4}\cdot n+\binom{n}{2}\right)\nonumber\\
&\ge\eta\binom{n}{2}-(\delta-\eta)\frac{\eta n^2}{4}\ge\frac{\eta}{2}\binom{n}{2}>\binom{\frac{\eta n}{4}}{2}=\binom{|G_i|}{2}\nonumber,
\end{align}
a contradiction. The proof is completed.
\eproof

\bproof[Proof of Proposition~\ref{RTprop}]
Let $\delta:=\varrho_\ell(r)$. Then for fixed $\alpha,\varepsilon>0$, there exist $n_0\in\mathbb{N}$ and a sequence of graphs $G_{n_0},G_{n_0+1},\dots$ such that for every positive integer $n\ge n_0$, $G_n$ is an $n$-vertex $K_r$-free graph with $e(G)\ge(\delta-\varepsilon)\binom{n}{2}$ and $\alpha(G)\le\alpha n$. For each $G_n$, we apply Lemma~\ref{RT2} with $\eta=\sqrt{\alpha}$ to obtain an $n'$-vertex subgraph $G_n'$ with $\delta(G_n')\ge(\delta-\varepsilon-\sqrt{\alpha})n'$ and $\alpha(G_n')\le\alpha n\le4\sqrt{\alpha}n'$, where $n'\ge\frac{\sqrt{\alpha}}{4}$. Then apply Lemma~\ref{RT1} to $G_n'$ to obtain a blow-up graph $G^*_n$ such that $|G^*_n|=n$, $\delta(G^*_n)\ge(\delta-2\varepsilon-\sqrt{\alpha})n$ and $\alpha(G^*_n)\le8\sqrt{\alpha}n$.

This sequence shows that $\varrho_\ell(r)-2\varepsilon\le \varrho^*_\ell(r)$. As this is for arbitrary $\varepsilon$, it shows that $\varrho_\ell(r)\le \varrho^*_\ell(r)$. Note that $\varrho_\ell(r)\ge\varrho^*_\ell(r)$ is trivial, and then we are done.
\eproof

\section{Concluding remarks}\label{sec7}

In our result, we determine the asymptotically optimal minimum degree condition for Problem~\ref{prob1.5} in the case $\ell\ge\frac{3r}{4}$, which is due to our limited understanding on the almost cover step (Lemma~\ref{tiling lem}). Note that by using the general lower bound of $\varrho^*_\ell(r)$ in Theorem~\ref{RTthm} and revisiting the calculation in Lemma~\ref{tiling lem}, we can show Lemma~\ref{tiling lem} for $\ell\ge0.72r$. Comparing with previous results, our result is firstly nontrivial for $r=8$ and $\ell=6$. On the other hand, Lemma~\ref{reachable lem} holds for all $\ell\ge 2$, which hopefully will be useful in the study of the remaining cases of Problem~\ref{prob1.5}.

\bibliographystyle{abbrv}
\bibliography{ref}

\end{document}